\documentclass[12pt]{article} \usepackage{amsmath,amssymb}

\begin{document} \def\square#1{\vbox{\hrule \hbox{\vrule\hbox to #1 
pt{\hfill}\vbox{\vskip #1 pt}\vrule}\hrule}} 
\newtheorem{guess}{Proposition }[section] 
\newtheorem{theorem}[guess]{Theorem} \newtheorem{lemma}[guess]{Lemma} 
\newtheorem{corollary}[guess]{Corollary} \def\Ext{\mathrm {Ext}} 
\def\Hom{\mathrm {Hom}} \def\Jac{\mathrm {Jac}} \def\deg{\mathrm {deg}} 
\def\dim{\mathrm {dim}} \def\End{\mathrm {End}} \def\cod{\mathrm 
{cod}}\def\C{\mathbb C}

\centerline {\Large \em The Atiyah--Jones conjecture
                         for rational surfaces}

\vspace{5mm} \centerline{\large \em Elizabeth Gasparim}

\begin{quote}{\small We show that if the Atiyah--Jones conjecture holds 
for a surface $X,$ then it also holds for the
 blow-up of $X$ at a point. 
Since the conjecture is known to hold for ${\mathbb P}^2$ and for ruled 
surfaces, it follows that the conjecture is true for all rational 
surfaces.} \end{quote}

  If $P \rightarrow X$ is a principal $SU(2)$ bundle over a Riemannian 
four-manifold $X$, with $c_2(P) = k >0$, and $A$ is a connection on $P$, 
the Yang-Mills functional $$YM(A)~=~\int_X||F_A||^2$$ is minimal precisely 
when the curvature $F_A$ is anti-self dual, i.e. $F_A = - *F_A$, in which 
case $A$ is called an instanton of charge $k$ on $X.$

Let ${\mathcal MI}_k(X)$ denote the  moduli space of framed instantons on $X$ with 
charge $k$ and let ${\mathcal C}_{k}(X)$ denote the space of all framed  gauge 
equivalence classes of connections on $X$ with charge $k.$ In 1978, Atiyah 
and Jones [AJ] conjectured that the inclusion ${\mathcal MI}_k(X) \rightarrow 
{\mathcal C}_{k}(X)$ induces an isomorphism in homology and homotopy through a 
range that grows with $k.$ The original statement of the
 conjecture was for the case when $X$ is a sphere, but the question 
readily generalises for other 4-manifolds.

The stable topology of these moduli spaces was understood in 1984, when 
Taubes [Ta] constructed instanton patching maps
 $t_k\colon {\mathcal MI}_k (X) \rightarrow {\mathcal MI}_{k+1} (X)$ and showed 
that the stable limit $\displaystyle\lim_{k\rightarrow \infty } {\mathcal 
MI}_k$ indeed has the homotopy type of $ {\mathcal C}_{k}(X).$ However, 
understanding the behaviour of the maps $t_k$ at finite stages is a finer 
question. Using Taubes' results, to prove the Atiyah--Jones conjecture it 
then suffices to show that the maps $t_k$ induce isomorphism in homology 
and homotopy through a range.
 
In 1993, Boyer, Hurtubise, Milgram and Mann [BHMM] proved that the 
Atiyah--Jones conjecture holds for the sphere $S^4$ and in 1995, Hurtubise 
and Mann [HM] proved that the conjecture is true for ruled surfaces.

In this paper I show that if the Atiyah--Jones conjecture holds true for a 
complex surface $X$ then it also holds for the surface $\widetilde{X}$ 
obtained by blowing-up $X$ at a point. In particular, it follows that the 
conjecture holds true for all rational surfaces.

Kobayashi and Hitchin gave a one-to-one correspondence between instantons 
on a topological bundle $E$ over $X$ and holomorphic structures on $E,$ 
see [LT]. Using this correspondence I translate the Atiyah--Jones 
conjecture into the language of holomorphic bundles and compare the 
topologies of the moduli spaces ${\mathfrak M}_k(X)$ and ${\mathfrak 
M}_k(\widetilde{X})$ of stable holomorphic bundles on $X,$ resp. 
$\widetilde{X},$ having $c_1=0$ and $c_2=k.$ More precisely, I will 
consider slightly enlarged moduli spaces ${\mathfrak 
M}^f_k(\widetilde{X})$ of bundles on $\widetilde{X}$ framed near, but not 
{\em on}, the exceptional divisor $\ell.$ That is, we fix a small tubular 
neighborhood $N(\ell)$ of $\ell$ in $\widetilde{X}$ and assign frames on 
$N^0:=N(\ell) - \ell.$ Existence of such frames is guarantied by 
[G2,Thm.\thinspace 4.1]. These
 framed moduli spaces ${\mathfrak M}^f_k(\widetilde{X})$ that are 
naturally related to the moduli spaces ${\mathfrak M}^f_k(X)$
 by a holomorphic gluing construction (Proposition 
4.1).

The structure of the proof is the following. First we give a concrete 
description of instantons on $\widetilde{{\mathbb C}^2}$ (section 1) and 
compute its numerical invariants (section 2). We then study moduli spaces 
${\mathfrak M}_k(\widetilde{X})$ of instantons or equivalently stable 
bundles on $\widetilde{X},$ with respect to a polarisation 
$\widetilde{\mathcal L} = N{\mathcal L}-{\ell}$ where ${\mathcal L}$ is a polarisation 
on $X.$ We show that removing the singular points of this moduli space 
does not affect homology in a range sufficient to our calculations and 
thereafter work only with its smooth points. Moreover, the direct image of 
a stable bundle might yield an unstable bundle, which does not correspond 
to an instanton on $X.$ We show that removing such unstable bundles from 
${\mathfrak M}_k({X})$ does not affect our homology calculations either 
(section 3).

We  show that any framed instanton on $\widetilde{X}$ is uniquely 
determined by
 holomorphic patching of a framed instanton on $X$ and an instanton on 
$\widetilde{{\mathbb C}^2}$ (section 4). We define 
framed moduli spaces (section 5) and then prove that the local moduli 
space ${\mathcal N}_i^f$ of (framed) instantons on
 $\widetilde{{\mathbb C}^2}$ with charge $i$ is a smooth complex manifold 
(section 6). As a consequence we have stratifications $${\mathfrak 
M}^f_k(\widetilde{X}) \simeq \bigcup_{i=0}^{k} {\mathfrak M}^f_{k-i}(X) 
\times {\mathcal N}^f_i.$$ We prove that the  map ${\mathfrak 
M}^f_k(\widetilde{X}) \rightarrow {\mathfrak M}^f_{k+1}(\widetilde{X})$  
obtained by translating  Taubes' map to bundles via Kobayashi--Hitchin 
correspondence 
is homotopic equivalent to a map that preserves the stratifications 
(section 7). Using Leray spectral sequences, 
we then show that
 Atiyah--Jones conjecture for $\widetilde{X}$ follows from the 
corresponding statement for $X.$

Acknowledgements: 
I thank R. J. Milgram and Tony Pantev for  enlightening
discussions, and J.P. Santos for pointing out a gap on an earlier 
version of this paper. 
 \vspace{3mm}
 
\noindent {\bf 1. Moduli of bundles on $\widetilde{{\C}^2}$ with fixed 
splitting type}

Let $\pi\colon \widetilde{X}\rightarrow X$ be the blow-up of a point $x$ 
on a compact
 surface $X$ and let $\ell $ be the exceptional divisor. In this paper all 
bundles have rank 2 and $c_1=0.$ Given a bundle $E$ over $X$ we pull it 
back to $\widetilde{X}$ and then modify it by gluing in new data near the 
exceptional divisor to construct a bundle $\widetilde{E}$ satisfying 
$\widetilde{E}\vert_{\widetilde{X}-\ell}\simeq E\vert_{X-x}.$ The 
difference of Chern classes $c_2(\widetilde{E}) -c_2(E)$ depends only on 
the restriction of $\widetilde{E}$ to a small neighborhood $N(\ell)$ of 
$\ell.$ To begin with, suppose that $N(\ell)$ is isomorphic to the blow-up 
of
 ${\C}^2$ at the origin, denoted $\widetilde{{\C}^2}.$ Given a rank two 
bundle $V$ on $\widetilde{{\C}^2}$ with vanishing first Chern class, there 
exists an integer $j \geq 0$ determined by the restriction of $V$ to the 
exceptional divisor, such that $V\vert_{\ell} \simeq {\mathcal O}_{{\mathbb 
P}^1}(j) \oplus {\mathcal O}_{{\mathbb P}^1}(-j).$ This $j$ is called the {\em 
splitting type} of the bundle. We denote by 
${\mathcal O}_{\widetilde{{\mathbb C}^2}}(j)$  the unique line bundle 
on $\widetilde{{\mathbb C}^2}$ having first Chern class $j,$ that is, 
the pull back of ${\mathcal O}_{{\mathbb P}^1}(j).$  
 We choose {\em canonical} coordinates for 
$N(\ell)= U \cup V,$ where $U \simeq {\C}^2 $ has coordinates $(z,u)$ and 
$V \simeq {\C}^2$ has coordinates $(\xi, v)$ with $(\xi,v)= (z^{-1},zu)$ 
in $U \cap V \simeq ({\C}-0) \times {\C}.$

\vspace{3mm}

\noindent {\sc Theorem} [G3] {\it Every holomorphic rank two bundle $V$ 
over $\widetilde{{\C}^2}$ with vanishing first Chern class and
 splitting type $j$ is an algebraic extension of the form $$0 \rightarrow 
{\mathcal O}_{\widetilde{{\C}^2}}(-j) \rightarrow V \rightarrow {\mathcal 
O}_{\widetilde{{\C}^2}}(j) \rightarrow 0$$ and the extension class can be 
represented in canonical coordinates by a polynomial of 
the form 
$$p= \sum_{i=1}^{2j-2}\sum_{l=i-j+1}^{j-1}p_{il}z^lu^i.$$ } \vspace{3mm}

It follows that the local (= over $N(\ell)$) moduli problem can be studied 
by considering extensions of line bundles modulo bundle isomorphism, even 
though this is definitely not the case for bundles on the compact surface 
$\widetilde{X}$. The fact that $V$ is determined by an extension of degree 
$2j-2$ implies that $V$ is determined over $N(\ell) $ by its restriction 
to the $2j-2$nd formal neighborhood of the exceptional divisor.
 Hence, our  assumption that $N(\ell) \simeq \widetilde{{\C}^2}$ is 
not restrictive. Consider the quotient space
 $${\mathcal M}_j\colon =
 \{ bundles \,\, on \,\, \widetilde{{\C}^2} \,\, with \,\, splitting \,\, 
type \,\, j \} / \sim,$$ where $\sim $ denotes bundle isomorphism. ${\mathcal 
M}_j$ is topologized as follows: 
In canonical coordinates, 
the extension class is expressed by the complex polynomial $p$ 
with $n=2j(j+1)$ 
coefficients  as in Theorem [G3]. In such canonical coordinates, 
elements of ${\mathcal 
M}_j$ are represented by transition matrices of the form
 $$\left( \begin{matrix} z^j & p \cr 0 & z^{-j} \end{matrix} 
\right).\eqno{(1)}$$
 We set
 $p\sim p'$ if the corresponding bundles are holomorphically equivalent. 
This construction identifies ${\mathcal M}_j $ to a quotient of ${\mathbb 
C}^n$ under a group action and endows ${\mathcal M}_j$ with the structure of a 
quotient stack. For more on the structure of ${\mathcal M}_j$ see [G3] and 
[G5]. The stack ${\mathcal M}_j$ has a corse moduli space
that can be decomposed into Hausdorff strata using the instanton
numerical invariants, see [BG2]. 
To study the moduli problem we calculate numerical
 invariants. \vspace{3mm}

\noindent {\bf 2. Computation of local numerical invariants}

In \S1 we saw that an element $V \in {\mathcal M}_j$ is determined by its 
splitting type $j$ and extension class $p.$ We write $V = V(j,p)$.  It is 
important to note that the bundle  $V$ is trivializable on the 
complement of $\ell.$ \vspace{3mm}

\noindent {\sc Lemma} [G2] {\it Every holomorphic bundle on 
$\widetilde{{\C}^2}$ with vanishing first Chern class is trivial on 
$\widetilde{{\C}^2}$ minus the exceptional divisor. } \vspace{3mm}
 
\noindent It follows that the local bundle $V(j,p)$ can be glued to any 
bundle $E$ pulled back from $X$ to form a new bundle $\widetilde{E}$ on 
$\widetilde{X}.$ We also note that the topology of $\widetilde{E}$ does 
not depend on the attaching map. In fact, the attaching is given by a 
holomorphic map $\psi \colon {\C}^2 - \{0\} \rightarrow Sl(2,\C).$ By 
Hartog's theorem $\psi$ extends to the origin, and therefore is  
homotopic to the identity;
 thus not contributing to the Chern numbers of $\widetilde{E}.$ Hence 
2 sets of data $(E,j,p, \psi)$ and $(E,j,p,\psi')$ determine the same 
topological bundle, although in general they determine distinct 
holomorphic bundles.

\vspace{3mm}

\noindent {\sc Lemma} [G4] {\it Every holomorphic rank two bundle 
$\widetilde{E}$ on $\widetilde{X}$ with $c_1(\widetilde{E})=0$ is 
topologically determined by a triple $(E,j,p)$ where $E$ is a bundle over 
$X,$ $j$ is a non-negative integer and $p$ is a polynomial.} \vspace{3mm}

\noindent Here $E= \pi_*(\widetilde{E})^{\vee\vee}$ and the Chern class 
difference $c= c_2(\widetilde{E})-c_2(E)$
 depends only on the local data $V(j,p).$ 
\vspace{3mm}

\noindent {\sc Definition 2.1} We call $c$ the {\it (local) charge} of 
the bundle $\widetilde{E}$ near the exceptional divisor. When 
considering the local situation, we also refer to $c$ as the 
{\it (local) charge} of $V.$ 
\vspace{3mm}

The terminology comes from 
instanton moduli spaces, where $c$ represents the change in topological 
charge of the instanton on $\widetilde{X}$ given by the patching of 
$V(j,p).$ Friedman and Morgan [FM] gave the bounds $j \leq 
c_2(\widetilde{E})-c_2(E) \leq j^2,$ and we proved 
sharpness. \vspace{3mm}

\noindent {\sc Theorem} [G1] {\it The bounds $j \leq 
c_2(\widetilde{E})-c_2(E) \leq j^2$ are sharp. \vspace{3mm}}

\noindent In [BG1] we used elementary transformations to show that
 all intermediate values occur. The following result implies non-emptiness 
of the strata appearing in the stratification in \S7. \vspace{3mm}

\noindent {\sc Theorem} [BG1] {\it For every integer $k$ satisfying $j 
\leq k \leq j^2$ there is a (semistable) 
holomorphic bundle on $\widetilde{X}$ with 
splitting type $j$ and such that $ c_2(\widetilde{E})-c_2(E) =k.$ 
\vspace{3mm}}

\noindent In section 7 we use two finer numerical invariants, defined as 
follows.
 An application of Riemann--Roch ( [FM], p.\thinspace 366)  gives 
$$c=c_2(\widetilde{E})-c_2(E)= l(R^1\pi_*\widetilde{E}) + l(Q),\eqno(2)$$ 
where $Q$ is the skyscraper sheaf defined by the exact sequence $$0 
\rightarrow \pi_*(\widetilde{E}) \rightarrow 
\pi_*(\widetilde{E})^{\vee\vee} \rightarrow Q \rightarrow 0,$$
and $l$ denotes length. The pair of 
local analytic invariants $l(Q)$ and $l(R^1\pi_*\widetilde{E})$ gives 
strictly finer information than the local charge $c,$ and has interesting 
properties such as:

\vspace{3mm}

\noindent {\sc Theorem} [BG2]{\it The pair of numerical invariants $l(Q)$ 
and $l(R^1\pi_*\widetilde{E})$ gives the coarsest stratification of ${\mathcal 
M}_j$ into maximal Hausdorff components. } \vspace{3mm}

\noindent
 Sharp bounds and nonemptiness of intermediate strata are given by the 
following results. \vspace{3mm}

\noindent {\sc Theorem}[G2]{\it Let $j>0$ be the splitting type of $E,$ 
then the following bounds are sharp
 $$j-1 \leq l (R^1\pi_* (E)) \leq j(j-1)/2 $$ and $$1 \leq l(Q)\leq 
j(j+1)/2.$$ The upper bounds occur only at the split bundle. If $j=0$ then 
both invariants are zero.} \vspace{3mm}

\noindent {\sc Theorem}[BG1]{\it For every pair of integers $(w,h)$
 satisfying
 $j-1 \leq h \leq j(j-1)/2 $ and $1 \leq w \leq j(j+1)/2$ with $j \geq 0$ 
there exists a rank 2 vector bundle $E$ on $\widetilde{{\C}^2} $ with 
splitting type $j$ having numerical invariants $l (R^1\pi_* (E)) = h$ and 
$l(Q) = w.$} \vspace{3mm}

Non-emptiness of intermediate strata is needed for the proof of Lemma 7.3. 
The remainder of this section gives some examples of how these invariants are 
calculated. Other examples of calculations of these invariants appear in 
[G1] and [BG2]. In [GS] a Macaulay2 program that calculates both 
invariants is given, and the following simple formula for 
$l(R^1\pi_*\widetilde{E})$ is proven.

\vspace{3mm} \noindent {\sc Theorem} [GS] {\it Let $m$ denote the largest 
power of $u$ dividing $p,$ and suppose $m>0.$ If $\widetilde{E}$ 
is the bundle defined 
by data $(E,j,p,\phi),$ then $$l(R^1\pi_*\widetilde{E})= \left(\begin{matrix} j 
\cr 2 \end{matrix}\right) - \left(\begin{matrix} j-m \cr 2 
\end{matrix}\right) .$$ }

\noindent The following example illustrates a computation of $l(Q).$ 
\vspace{3mm}

\noindent {\sc Example 2.2 } {Let $E$ be given by data $(j,z^nu); 
n\geq 1.$ We show that $l(Q) = {\tiny \left ( \begin{matrix} n \cr 2 
\end{matrix} \right)}, $ it then follows from theorem [GS] and equality 
(2) that the charge is $$c(E)= l(Q)+ l(R^1\pi_*\widetilde{E})= \left( 
\begin{matrix} n \cr 2 \end{matrix} \right) + \left(\begin{matrix} j \cr 2 
\end{matrix}\right) - \left(\begin{matrix} j-1 \cr 2 \end{matrix}\right)=
 \left( \begin{matrix} n \cr 2 \end{matrix} \right) +j. $$

We use the method of [BG2] section 6.2. Let $M = (\pi_* E)_0^{\wedge}$ 
denote the completion of the stalk
 $(\pi_*E)_0.$ Let $\rho$ denote the natural inclusion of $M$ into its 
double dual $\rho: M \hookrightarrow M^{\vee\vee}.$ We want to compute $l(Q) = 
dim\, coker(\rho).$ By the theorem on formal functions $$M \simeq 
\lim_{\longleftarrow}H^0(\ell_n, \widetilde{V}|{\ell_n}),$$
where $\ell_n$ denotes the $n-$th infinitesimal neighborhood of $\ell.$
 There are 
simplifications that make it easy to calculate $M,$ (cf. [BG2] or [G1]). 
So,
 to determine $M$ it suffices to
 calculate $H^0(\ell_{2j-2}, \widetilde{V}|{\ell_{2j-2}}),$
 and the relations among its generators under the action of ${\mathcal 
O}^{\wedge}_0 ( \simeq {\bf C}[[x,y]]$). In this example, $E$ is given by 
transition matrix ${\tiny \left(\begin{matrix} z^j & z^nu \cr 0 & z^j 
\end{matrix} \right).}$ We find that $M = {\bf C}[[x,y]]
\langle\beta_0,\beta_1, 
\cdots, \beta_{n}, \gamma\rangle $ where, for $1 \leq i \leq n,$

$$ \beta_i = \left(\begin{matrix} -z^iu \cr 
z^{j-n+i}\end{matrix}\right),\,\,\, \gamma = \left(\begin{matrix} 0 \cr 
u^{n-1} \end{matrix}\right) $$ with relations $$ x\beta_{i} - y 
\beta_{i-1} =0.$$ Consequently, $M^{\vee} = \langle B , C \rangle$ 
is free on the  
generators $$ B :\left\{\begin{array} {l} \beta_i \rightarrow x^iy^{n-i} 
\cr \gamma \rightarrow 0 \cr \end{array}\right.\,\,\,
 C :\left\{\begin{matrix} \beta_i \rightarrow 0 \cr \gamma \rightarrow 1 
\cr \end{matrix}\right..
 $$
 $M^{\vee \vee} = \langle {\mathcal B}, {\mathcal C} \rangle $
 is free on the generators $$ {\mathcal 
B} :\left\{\begin{matrix} B \rightarrow 1 \cr C \rightarrow 0 \cr 
\end{matrix}\right.\,\,\, \,\,\,
 {\mathcal C} :\left\{\begin{matrix} B \rightarrow 0 \cr C \rightarrow 1 \cr 
\end{matrix}\right..$$ The inclusion into the bidual $\rho :M \rightarrow 
M^{\vee \vee},$ takes $x$ into evaluation at $x,$ therefore $$ \rho 
:\left\{\begin{array}{l} \beta_i \rightarrow x^iy^{n-i} {\mathcal B} \cr 
\gamma \rightarrow {\mathcal C} \end{array}\right.$$
 and $coker( \rho) = \langle x^iy^{n-i-1}{\mathcal B}; 0 \leq i \leq n \rangle.$ We 
conclude that $$l(Q) = dim\, coker( \rho) = {n \choose 2}.$$} \vspace{1mm}

\noindent {\bf 3. Moduli of bundles on $\widetilde{X}$}

For the moduli problem on $\widetilde{X}$ we need stability conditions.
 If ${\mathcal L}$ is an ample divisor on $X$ then, for large $N,$ the divisor 
$\widetilde{\mathcal L} = N{\mathcal L}-{\ell}$ is ample on $\widetilde{X}.$ We 
fix, once and for all,
 polarisations ${\mathcal L}$ and $\widetilde{\mathcal L}$ on $X$ and 
$\widetilde{X}$ respectively.
 For a polarised surface $Y,$ the notation
 ${\mathfrak M}_k(Y)$ stands for moduli of rank 2
 bundles on $Y$ having $c_1=0,$ $c_2=k$ 
and  slope stable with respect to the fixed polarisation.
 If $E$ is ${\mathcal L}$-stable
 on $X,$ then $\pi^*(E)$ is $\widetilde{\mathcal L}$-stable on $\widetilde{X}$ 
[FM, Thm.\thinspace 6.5]. Hence, the pull back map induces an inclusion of 
moduli spaces
 ${\mathfrak M}_k(X) \hookrightarrow {\mathfrak M}_k(\widetilde{X}).$
 
Our objective is to show that given a map ${\mathfrak M}_k(X) \rightarrow 
{\mathfrak M}_{k+1}({X})$ inducing a homology equivalence through a range, 
there is a map ${\mathfrak M}_k(\widetilde{X}) \rightarrow {\mathfrak 
M}_{k+1}(\widetilde{X})$ inducing a homology equivalence through a 
comparable range, and that both ranges go to $\infty$ with $k.$ Given that 
the homology equivalence for the case of $X={\mathbb P}^2$ holds up to 
dimension $\lfloor k/2 \rfloor-1,$ we expect to obtain an equivalence 
range at best equal to this one. From here on, we are free to modify 
${\mathfrak M}_k(X)$ and ${\mathfrak
  M}_k(\widetilde{X})$ in any way that does not change topology up to 
dimension $k/2=c_2/2$. Firstly, we remove singularities of the moduli 
spaces; secondly, we remove semistable bundles on ${\mathfrak M}_i(X)$ 
that are not stable.

\vspace{3mm}

\noindent{ \em 3.1 Removing singularities}

We remove the singular points of the moduli spaces, in order to work only 
with smooth manifolds. We use the following results of Kirwan and 
Donaldson. \vspace{3mm}

\noindent {\sc Theorem} ([Ki],\thinspace Cor. 5.4) {\it Let $X$ be a 
quasi-projective variety and $m$ a non-negative integer such that every 
$x_0 \in X$ has a neighborhood in $X$ isomorphic to $$\{x \in {\C}^N | 
f_1(x)= \cdots = f_M(x)=0\}$$ for some integers $N,$ $M,$ and holomorphic 
functions $f_i$ depending on $x_0$ with $M \leq m.$ If $Y$ is a closed 
subvariety of codimension $k$ in $X,$ then for $ q < k-m,$} 
$$H_q(X-Y,{\mathbb Z}) \simeq H_q(X,{\mathbb Z}).$$ Given a complex 
surface $S$ with polarisation $H,$ let $\Sigma_k \subset {\mathfrak 
M}_k(S)$ denote the algebraic subvariety representing
 bundles $E$ with $H^2(\End_0 E) \neq 0.$ A local model of the moduli 
space is determined by the kernel of the Kuranishi map $Sym^2 H^1(\End_0 E) 
\rightarrow H^2(\End_0 E)$ parameterising small deformations of $E.$ In 
case $H^2(\End_0 E)=0,$ small deformations of $E$ are unobstructed and $E$ 
is a smooth point of the moduli space. Therefore, singular points satisfy 
$H^2(\End_0 E) \neq 0$ and the singular part of ${\mathfrak M}_k$ is 
contained in $\Sigma_k.$

\vspace{3mm}

\noindent {\sc Theorem} ([Do],\thinspace Thm. 5.8) {\it There are 
constants $a,b$ depending only on $S$ and the ray spanned by $H$ in 
$H_2(S)$ such that:} $$\dim_{\mathbb C}\Sigma_k \leq a+b\sqrt{k}+3k.$$ 
\noindent {\sc Proposition 3.1} {\it Removing the 
singularity set $Sing$ of ${\mathfrak M}_k(S)$ does not change homology in 
dimension less than k. That is, for $q<k$} $$H_q({\mathfrak M}_k(S))=H_q 
({\mathfrak M}_k(S)-Sing).$$

\noindent {\it Proof.} By Kuranishi theory, points
 $E \in {\mathfrak M}_k(S)$ satisfying $H^2(\End_0 E) = 0$ are smooth 
points. Therefore, the singularity set of $ {\mathfrak M}_k(S)$ is 
contained in $\Sigma_k.$ Moreover, the moduli space is defined on a 
neighborhood of a singular point by at most $\dim H^2(\End_0 E)$ equations. In 
fact, in the proof of his theorem cited above, Donaldson shows that $ \dim 
H^2(\End_0 E) \leq a+b\sqrt{k}+3k.$ Therefore, using Kirwan's result we 
have that $H_q({\mathfrak M}_k(S))=H_q ({\mathfrak M}_k(S)-Sing)$ for 
$q<\dim {\mathfrak M}_k(S) - 2(a+b\sqrt{k}+3k) < k.$ 
 \hfill\square{5} \vspace{3mm}

 Henceforth we consider only the smooth part of ${\mathfrak 
M}_k(\widetilde{X}),$ which in what follows stands for just the set of its 
smooth points. Similarly, when considering ${\mathfrak M}_k(X)$ we will 
consider only its set of smooth points, keeping in mind that this does not 
change the homology up to dimension $k.$ \vspace{3mm}

\noindent {\em 3.2 Removing unstable bundles}

Given an $\widetilde{\mathcal L}$-stable bundle $\widetilde{E},$ the bundle
 $E=(\pi_* \widetilde{E})^{\vee\vee}$ is ${\mathcal L}$-semistable [FM, 
Thm.\thinspace 6.5]. We want to remove the set of strictly semistable 
bundles that appear in this process, leaving only stable bundles over $X,$ 
as these are the ones corresponding to instantons. \vspace{3mm}

\noindent {\sc Proposition 3.2} {\it Removing  
semistable bundles that are not 
stable does not change homology in dimension less than $7k-c.$ That is 
if  ${\mathfrak M}^{su}$ denotes the subset of semistable bundles on
 $X$ that are not stable, then there is a constant $c$
depending only on $X$ such that for $q< 7k -c$
$$H_q({\mathfrak M}_k(X))=H_q 
({\mathfrak M}_k(X)-{\mathfrak M}^{su}_k(X)).$$
} 

\noindent {\it Proof.} Let $E$ be a semistable bundle on $X$ that is not 
stable. Since $c_1(E)=0$
 there is a destabilising (saturated) 
line bundle $L$ with $\deg(L)=0.$ Hence, $E$ is 
an extension $$0\rightarrow L \rightarrow E \rightarrow {\mathcal F} 
\rightarrow 0,\eqno(3)$$ where ${\mathcal F}$ is a rank 1 
torsion free sheaf satisfying 
${\mathcal F}^{\vee\vee}=L^{-1}$ and fitting into a short exact sequence 
$$0\rightarrow {\mathcal F} \rightarrow {\mathcal F}^{\vee\vee} \rightarrow Q 
\rightarrow 0,$$ with $Q$ supported only at points.
Taking $\Hom(.,L)$ on the short exact sequence $(3)$ 
gives 
%$$0 \rightarrow  \Hom(Q,L) \rightarrow \Hom({\mathcal F}^{\vee\vee},L) 
%\rightarrow \Hom({\mathcal F},L) \rightarrow$$ 
$$\rightarrow \Ext^1(Q,L)
\rightarrow \Ext^1({\mathcal F}^{\vee\vee},L) \rightarrow 
\Ext^1({\mathcal F},L) \rightarrow \Ext^2(Q,L) \rightarrow
\Ext^2({\mathcal F}^{\vee\vee}, L) \rightarrow.$$ 
But, because $Q$ is supported at $k = c_2(E)$ points, we 
have $\Hom(Q,L)=\Ext^1(Q,L)=0,$ and $\Ext^2(Q,L) = {\mathbb C}^k,$ and 
consequently the long exact sequence becomes
 
$$0 \rightarrow H^1(L^2) \rightarrow \Ext^1({\mathcal F},L) 
\stackrel{\delta} {\rightarrow}
 Ext^2(Q,L) ={\mathbb C}^k \rightarrow H^2(L^2) \rightarrow .\eqno(4)$$
Note that not all sequences of the form $(3)$ 
result in a bundle $E,$ but as we only consider bundles, we must identify 
when the middle term is locally free. 

Locally we have a resolution
$$0 \rightarrow L \rightarrow E \rightarrow F^{\vee\vee} \rightarrow
Q \rightarrow 0$$
$$0 \rightarrow {\mathcal O} \rightarrow  {\mathcal O}^2 \rightarrow  {\mathcal O}
\rightarrow Q \rightarrow 0$$
and it follows that  $Q$ is a local complete intersection.
We write
$\displaystyle Q= {\mathcal O}/(p,q)$ with $p,q \in {\mathcal O}$ coprime.
Assuming this, we have that ${\mathcal E}xt^1(Q,{\mathcal O}) = 0 $ and
${\mathcal E}xt^2(Q,{\mathcal O}) \simeq        Q.$
The extension $(3)$ 
%$$0 \rightarrow L \rightarrow E \rightarrow F \rightarrow 0$$
has a corresponding class 
$[E] \in {\mathcal E}xt^1(F,L) \simeq {\mathcal E}xt^2(Q,L) =
 {\mathcal E}xt^2(Q,{\mathcal O}) \otimes L .$
In order for $E$ to be a bundle we need that $[E]$
have a non-zero value at each point in the support of $Q.$ We have 
$$\mbox{Ext}^1(F,L) \rightarrow \Gamma {\mathcal E}xt^1(F,L) =
   \Gamma {\mathcal E}xt^2(Q,L)= \mbox{Ext}^2(Q,L)$$
and the map from leftmost to rightmost terms in the above
expression is just the connecting homomorphism $\delta$ 
coming from
$$0 \rightarrow F \rightarrow F^{\vee\vee} \rightarrow Q \rightarrow 0.$$

 It follows that the 
dimension of the set ${\mathfrak M}^{su} : = {\mathfrak M}^{ss} - 
{\mathfrak M}^s$ of semistable bundles on $X$ that are not stable is given 
by $$\dim {\mathfrak M}_k^{su}= \dim(\delta(\Ext^1({\mathcal F},L) )
 \times \Jac(X)) $$ and using the exact sequence $(4)$ we get  
$$\dim{\mathfrak M}_k^{su}\leq k+ \dim(H^1(L^2) 
 \times \Jac(X)).$$ 
%It is a basic fact that dim$\Jac(X)= h^1(X,{\mathcal O}_X)$ is the 
%irregularity of $X.$
where $ c=\dim(H^1(L^2) 
 \times \Jac(X))$ depends only on $X.$
Now apply Kirwan's theorem with codimension 
$\dim{\mathfrak M}_k^s - \dim{\mathfrak M}_k^{su}= 8k-(1-b_1+b^+) - (k+c)$
(singularities were removed in the previous lemma).
 \hfill\square{5} \vspace{3mm}

Based on the results of this section, we assume in what follows that
 the moduli spaces ${\mathfrak M}_i$
 are smooth and contain only stable points. We now proceed to the study of 
framed bundles. \vspace{3mm}

\noindent {\bf 4. Holomorphic instanton patching}

 In this section we give the detailed construction of holomorphic 
instanton patching which gives the following result. \vspace{3mm}

\noindent {\sc Proposition 4.1} {\it Every (framed) instanton on 
$\widetilde{X}$ is obtained by holomorphic patching an instanton on $X$ 
to an instanton on $\widetilde{{\C}^2}.$} {\vspace{3mm}

We give the patching in terms of holomorphic bundles, via the 
Kobayashi--Hitchin correspondence. For a K\"ahler surface $X,$ Kobayashi 
and Hitchin gave a one-to-one correspondence between irreducible $SU(2)$ 
instantons of charge $k$ on $X$ and stable rank 2 holomorphic bundles $E$ 
on $X$ with $c_1=0$ and $c_2=k,$ see [LT]. For the noncompact surface 
$X=\widetilde{{\C}^2}$ this correspondence takes an instanton to a 
holomorphic bundle with an added trivialisation at infinity [Kn]. We show 
that a framed bundle on $\widetilde{X}$ is uniquely determined by a pair 
of framed bundles on $X$ and $\widetilde{{\C}^2}.$ 
\vspace{3mm}

\noindent{\sc Remark:} Note that all gluing maps  used here 
for  holomorphic patching of framed bundles are homotopic to the 
identity, this is a consequence of Hartog's theorem as explained in section 2.
\vspace{3mm}

We fix a neighborhood 
$N(\ell)$ of the exceptional divisor inside of $\widetilde{X}$ and
 write $\widetilde{X}= (\widetilde{X}-\ell) \cup N(\ell).$ Set $$N^0 : 
= N(\ell) -\ell =(\widetilde{X}-\ell) \cap N(\ell) .$$ The blow-up map 
gives an isomorphism $i_1\colon \widetilde{X}-\ell
 \rightarrow X - \{x\}.$ Based on the results of \S1, we 
we know that moduli of bundles on $N(\ell)$ are isomorphic to 
moduli of bundle on $\widetilde{{\mathbb C}^2}.$ 
Hence, to simplify the exposition,  we can assume that 
there is an isomorphism $i_2\colon N(\ell) \rightarrow \widetilde{{\C}^2} 
.$ Over the intersection $N^0$ we have isomorphisms $ N(x)-\{x\} 
\stackrel{i_1}{\leftarrow} N^0 \stackrel{i_2}{\rightarrow}
 \widetilde{{\C}^2} - \ell ,$ which we still denote by the same 
letters $i_1$ and $i_2.$ Using this isomorphisms we write framed bundles 
on $\widetilde{X}$ as pairs of framed bundles on $X$ and $ 
\widetilde{{\C}^2}.$ \vspace{3mm}

\noindent {\sc Lemma 4.2} {\it There is a one-to-one correspondence 
between algebraic vector bundles on $\widetilde{X}-\ell$ and algebraic 
vector bundles on $X.$} \vspace{3mm}

\noindent{\it Proof.} It is trivial that there is a one-to-one 
correspondence between bundles on $\widetilde{X}-\ell$ and
 bundles on $X-\{x\},$ because these are isomorphic surfaces. If $E$ is an 
algebraic bundle on
 $X-\{x\}$ it extends uniquely to a bundle on $X$ as follows. One extends 
$E$ over $x$ as a coherent sheaf ${\mathcal E}$ and then taking double dual 
one obtains a reflexive sheaf ${\mathcal E}^{\vee\vee}$ over $X.$ Since 
reflexive sheaves have singularities in codimension 2, it follows that 
${\mathcal E}^{\vee\vee}$ is locally free. We have ${\mathcal 
E}^{\vee\vee}\vert_{X-\{x\}} \simeq {\mathcal E}\vert_{X-\{x\}} \simeq E,$ so
 ${\mathcal E}^{\vee\vee}$ gives the vector bundle extension of $E.$ 
\hfill\square{5} \vspace{3mm}

\noindent {\sc Definitions 4.3} \begin{itemize}

\item Let $\pi_F \colon F \rightarrow Z$ be a bundle over a surface $Z$ 
that is trivial over $Z_0:= Z - Y,$ where 
$Y$ is a closed submanifold of $Z.$ Given two pairs $f=(f_1,f_2)\colon Z_0 
\rightarrow \pi_{F}^{-1}(Z_0)$ and $g=(g_1,g_2)\colon Z_0 \rightarrow 
\pi_{F}^{-1}(Z_0)$ of fibrewise linearly independent holomorphic sections 
of $F\vert_{Z_0},$
 we say that $f$ is {\em equivalent} to $g$ (written $f \sim g$)
if $\phi:=g \circ f^{-1} \colon V\vert_{Z^0} \rightarrow V \vert_{Z^0}$ 
extends to a holomorphic map  $\phi\colon F\rightarrow F$
  over the entire $Z.$ A {\it 
frame} of $F$ over $Z_0$ is an equivalence class of fibrewise linearly 
independent holomorphic sections of $F$ over $Z_0.$ The set of 
such frames 
$$Fram(Z_0,F) := Hol (Z_0, SL(2,{\mathbb C}))/ \sim$$ 
carries the quotient topology.  

\item
 A {\em framed } bundle $\widetilde{E}^f$ on $\widetilde{X}$ is a pair 
consisting of a bundle $\pi_{\widetilde{E}}\colon \widetilde{E} 
\rightarrow \widetilde{X}$ together with a frame of $\widetilde{E}$ over 
$N^0:=N(\ell) - \ell.$

\item A {\em framed } bundle $V^f$ on $\widetilde{{\C}^2}$ is a pair 
consisting of a bundle $\pi_V\colon V \rightarrow \widetilde{{\C}^2}$ 
together with a frame of $V$ over $\widetilde{{\C}^2} - \ell.$

\item A {\em framed } bundle $E^f$ on $X$ is a pair consisting of a bundle 
$E \rightarrow X$ together with a frame of $E$ over $N(x)-\{x\},$ where 
$N(x)$ is a small disc neighborhood of $x.$ We will always consider $N(x)= 
\pi_{\widetilde{E}}(N(\ell)).$ \end{itemize}

\vspace{3mm}

\noindent{\sc Proposition 4.4} {\it An isomorphism class 
$[\widetilde{E^f}]$ of a framed bundle on $\widetilde{X}$ is uniquely 
determined by a pair of isomorphism classes of framed bundles $[E^f]$ on 
$X$ and $[V^f]$ on $\widetilde{{\C}^2}.$ We write 
$\widetilde{E}^f=(E^f,V^f).$} \vspace{3mm }

\noindent {\it Proof.} By construction
 $ \widetilde{E} = E_{\amalg_{(s_1,s_2)=(t_1,t_2)}} V$ is made by 
identifying the bundles as well as the sections over $N^0$, so that the 
bundles satisfy $\widetilde{E}\vert_{N^0}= i_1^*(E\vert_{N(x)-\{x\}})= 
i_2^*(V\vert_{\widetilde{{\C}^2}-\ell})$ and the framing $(f_1,f_2)$ 
of $\widetilde{E}$ satisfies $(f_1,f_2)= (s_1,s_2)\circ i_1= (t_1,t_2) 
\circ i_2.$

Let $\phi\colon E^f \rightarrow {E'}^f$ be an isomorphism such that 
$\phi\circ (s_1,s_2)=(s'_1,s'_2)$ and let $\xi\colon V^f \rightarrow 
{V'}^f$ be an isomorphism such that $\xi\circ (t_1,t_2)=(t'_1,t'_2).$ We 
have the following diagram of bundle maps:

 $$\begin{matrix} \widetilde{E^f}\vert_{\widetilde{X}-\ell} & 
\rightarrow &
 E^f & \stackrel{\phi}{\rightarrow} & {E'}^f & \leftarrow & 
\widetilde{{E'}^f}\vert_{\widetilde{X}-\ell} \cr
  \downarrow \pi_{\widetilde{E}} & & \downarrow\pi_E & & \downarrow 
\pi_{E'} & & \downarrow \pi_{\widetilde{E'}} \cr \widetilde{X}-\ell & 
\stackrel {i_1}{\rightarrow} & X-\{x\}& = & X-\{x\} & \stackrel 
{i_1}\leftarrow & \widetilde{X}-\ell\cr \end{matrix}$$ Hence, 
$$\widetilde{E'}^f\vert_{\widetilde{X}-\ell} = 
i_1^*({E'}^f\vert_{X-\{x\}})
 = i_1^*\circ \phi(E^f\vert_{X-\{x\}}) =
 i_1^*\circ \phi \circ i_{1*} ( 
\widetilde{E}^f\vert_{\widetilde{X}-\ell})\eqno{(4)}$$ showing that $ 
i_1^*\circ \phi \circ i_{1*}$ is an isomorphism of $\widetilde{E}$ and 
$\widetilde{E'}$ over $\widetilde{X}-\ell$ such that $$ \phi \circ 
(f_1,f_2)= \phi\circ (s_1,s_2)\circ i_1 = (s'_1,s'_2)\circ i_1 = 
(f'_1,f'_2).\eqno{(5)}$$ On the other hand we have the second diagram of 
bundle maps:

 $$\begin{matrix} \widetilde{E^f}\vert_{N(\ell)} & \rightarrow &
 V^f & \stackrel{\xi}{\rightarrow} & {V'}^f & \leftarrow & 
\widetilde{{E'}^f}\vert_{N(\ell)} \cr
  \downarrow \pi_{\widetilde{E}} & & \downarrow\pi_V & & \downarrow 
\pi_{V'} & & \downarrow \pi_{\widetilde{V'}} \cr N(\ell) & \stackrel 
{i_2}{\rightarrow} & \widetilde{{\C}^2}& = & \widetilde{{\C}^2}& \stackrel 
{i_2}\leftarrow & N(\ell)\cr \end{matrix}.$$ Therefore,
 $$\widetilde{E'}^f\vert_{N(\ell)} = 
i_2^*({V'}^f\vert_{\widetilde{{\C}^2}}) = i_2^*\circ \xi 
(V^f\vert_{\widetilde{{\C}^2}}) =
 i_2^*\circ \xi \circ i_{2*} ( 
\widetilde{E}^f\vert_{N(\ell)})\eqno{(6)},$$ showing that $i_2^*\circ \xi 
\circ i_{2*}$ is an isomorphism of $\widetilde{E}$ and $\widetilde{E'}$ 
over $N(\ell)$ such that $$ \xi \circ (f_1,f_2)= \xi \circ (t_1,t_2)\circ 
i_2= (t'_1,t'_2) \circ i_2= (f'_1,f'_2).\eqno{(7)}$$ These isomorphisms 
agree over the intersection $N^0,$ in fact, by (4) and (6) $$ i_1^*\circ 
\phi \circ i_{1*} (\widetilde{E^f}\vert_{N^0}) = i_1^*\circ \phi 
(E^f\vert_{N(x)})=i_1^* ({E'}^f\vert_{N(x)}) =i_2^* ( 
{V'}^f\vert_{\widetilde{{\C}^2}-\{0\}})=$$ $$ = i_2^*\circ \xi 
(V^f\vert_{\widetilde{{\C}^2}-\{0\}})= i_2^* \circ \xi \circ i_{2*} 
(\widetilde{E^f}\vert_{N^0})$$ and moreover they also preserve the 
framings over the intersection, since over $N^0,$ we have, by (5) and (7) 
$$ \phi \circ (f_1,f_2)=
 (f'_1,f'_2)= \xi \circ (f_1,f_2).$$ By the gluing lemma this gives an 
isomorphism over the entire $\widetilde{X}$ and we get $\widetilde{E'} 
\simeq \widetilde{E}.$ \hfill\square{5} \vspace{3mm}

\noindent Note: It is also possible to define framings as being 
just trivializing sections, without putting the equivalence relation
at first, and later divide by automorphisms of  the bundle that 
preserve the framings. However, we found it convenient not to carry on too 
many inequivalent framings, so that the framed local moduli 
${\mathcal  N}_i^f$ remain 
finite dimensional.

\vspace{3mm}

\noindent {\bf 5. Framed moduli space of bundles and connections.}

Anti-self-dual connections on a connected sum have been extensively 
studied in the literature, cf. [DK,Chap.\thinspace 7] and [Ta2] for 
detailed expositions. In particular, results apply to connections on a 
blow-up, since differentiably $\widetilde{X} \sim X \# 
\overline{\mathbb{CP}^2}.$ A posteriori, applying the Kobayashi--Hitchin 
correspondence, one knows that over K\"ahler surfaces, the theory of 
stable bundles must be completely parallel to the theory of irreducible 
ASD connections. Adding frames rigidifies the problem and has the 
advantage of greatly simplifying gluing constructions. Note that a framing 
of a bundle induces a framing of the corresponding connection.

 To relate anti-self-dual connections on $X$ and $\widetilde{X}$ one needs 
some compatibility between the metrics $\widetilde{g}$
 on $\widetilde{X}$ and $g$ on $X.$ Since anti-self-duality is a 
conformally invariant condition, it suffices to give a conformal metric 
$\widetilde{g}$ on $\widetilde{X}$ such that 
$\widetilde{g}\vert_{\widetilde{X}- \ell}$ is conformally equivalent to 
$g\vert_{X-\{x\}}.$ Construction of such conformal metrics is carried out 
in detail in [TA2, p.\thinspace 65]. In what follows we consider
 fixed such metrics.  Our proof of the Atiyah--Jones conjecture works for 
all metrics $\widetilde{g}$ on $\widetilde{X}$ for which the conjecture 
holds true for the corresponding $g$ on $X.$
 
We now recall some well known facts about moduli of vector bundles. Given 
a $C^{\infty}$ complex vector bundle $F$ on a compact complex surface $Z,$ 
a holomorphic structure on $F$ is a $\bar\partial$ operator of type 
$(0,1)$ which is integrable, hence a holomorphic bundle on $Z$ is a pair 
${\mathcal F} = (F,\bar\partial );$ we denote the set of all these by ${\mathcal 
H}.$ Let ${\mathcal G}$ be the group of $C^{\infty}$-automorphisms of $F.$ 
Then ${\mathcal G}$ acts on ${\mathcal H},$ and $\bar\partial_1$ and 
$\bar\partial_2$ are in the same orbit of ${\mathcal G}$ if and only if the 
corresponding holomorphic bundles $(F,\bar\partial_1)$ and 
$(F,\bar\partial_2)$ are isomorphic.  Given a polarisation $H$ on $Z$ let 
${\mathcal H}^s_k(Z):={\mathcal H}^s(F) $ bet the subset of
 $H$-stable bundles (having $c_1=0$ and $c_2=k$);
 the corresponding moduli space
 is the quotient ${\mathcal H}^s_k/{\mathcal G}.$ By Maruyama's theorem this 
moduli space is a quasi projective variety [M].

We fix compatible polarisations ${\mathcal L}$ and $\widetilde{\mathcal L} = 
N{\mathcal L}-{\ell}$ on $X$ and $\widetilde{X}$ throughout as in section 3. 
Note that this is equivalent to considering the projective embeddings 
determined by ${\mathcal L}$ and $\widetilde{\mathcal L} $ and then,  regarding $X$ 
and $\widetilde{X}$ with the corresponding induced K\"ahler metrics $g$ 
and $\widetilde{g}$ to use stability with respect to
 the polarisations given by the K\"ahler classes.  We denote the moduli 
spaces of ${\mathcal L}$ and $\widetilde{\mathcal L}$ stable bundles with $c_1=0$ 
and $c_2=k$ just by ${\mathfrak M}_k(X)$ and ${\mathfrak 
M}_k(\widetilde{X})$ respectively. Moreover, using the results of section 
3, we removed all singularities, and accordingly we assume that 
${\mathfrak M}_k(\widetilde{X})$ and ${\mathfrak M}_k(X)$ are smooth.

 Frames are added as in definition 4.3, and an isomorphism between framed 
bundles $(E_1,f^1)$ and $ (E_2,f^2)$ is an isomorphism $\phi \in {\mathcal 
G}\colon E_1 \rightarrow E_2 $ taking $f^1$ to $f^2,$ that is, such that 
$f^2 = f^1 \circ \phi.$ Hence, the complex gauge group acts on framed 
bundles and the framed moduli spaces are obtained as quotient of such 
actions. Explicitly, if ${\mathcal E}$ and $\widetilde{\mathcal E}$
 having $c_1=0$ and 
$c_2=k$ 
are the differentiable supports for our holomorphic bundles over $X$ and 
$\widetilde{X},$ with ${\mathcal G}(\mathcal E)$ and 
${\mathcal G}(\widetilde{\mathcal E})$ as their 
groups of automorphisms, then the framed moduli spaces are: \vspace{3mm}

\noindent {\sc Definitions 5.1} \begin{itemize}

\item ${\mathfrak M}_k^f(\widetilde{X}):= 
\left\{(\bar \partial_{\widetilde{\mathcal E}},f): 
f \in Fram\left ( N^0,(\widetilde{\mathcal E},
\bar \partial_{\widetilde{\mathcal E}} )\right ) \right\}/
{\mathcal G}(\widetilde{\mathcal E})$

\item ${\mathfrak M}_k^f(X):= 
\left\{(\bar \partial_{{\mathcal E}},f): 
f \in Fram\left ( N(x) - \{x\},({\mathcal E},
\bar \partial_{{\mathcal E}} )\right ) \right\}/
{\mathcal G}({\mathcal E})$

\end{itemize}

Note that ${\mathfrak M}^f_k(\widetilde{X})$ maps to ${\mathfrak 
M}_k(\widetilde{X}),$ and ${\mathfrak M}_k^f(X)$ maps to ${\mathfrak 
M}_k(X)$
 by projection onto the first coordinate. 
In section 6, we consider the local moduli spaces of bundles
defined only over a small 
tubular neighborhood $N(\ell)$ of the divisor, or what is equivalent, 
bundles on $\widetilde{{\mathbb C}^2}.$
 In 6.1 we explicitly calculate  dimensions of framed 
local moduli. This  clarifies 
the important point that the strata appearing in $(7)$ have 
finite codimension inside the global framed moduli spaces. 
An independent calculation of codimensions is given in 7.3 as well.
\vspace{3mm}

\noindent {\sc Remark 5.2} We remark that is it possible to 
choose framings in a different way, so that
moduli spaces of framed instantons 
remain finite dimensional. In fact, 
by Theorem [G3] bundles on a tubular neighborhood 
of the exceptional divisor are completely determined by 
their restriction to a finite 
infinitesimal neighborhood of order $2j-2.$ 
Consequently one could choose framings only on such finite 
infinitesimal neighborhoods and work with finite dimensional 
framed moduli spaces.

%Here again we consider only %smooth points, according to section 3.

\vspace{3mm}

%\noindent {\sc Remark 5.2}:
% ${\mathfrak M}^f_k(\widetilde{X}),$ % ${\mathfrak M}^f_k(X), {\mathcal M}^f$ 
%and ${\mathcal N}^f$ %denote the framed versions of 
% ${\mathfrak M}_k(\widetilde{X}),$ ${\mathfrak M}_k(X),
% {\mathcal M}$ % and ${\mathcal N}$ 
%respectively, the first 2 containing only %smooth points and stable 
%bundles as per section 3, %and all 4 spaces framed according to the 
%definitions %of section 4. \vspace{3mm}

\noindent {\bf 6. Local moduli}

For the moduli problem we need to consider the spaces ${\mathcal N}_i^f$ of 
isomorphism classes of framed bundles on the neighborhood 
$N(\ell)\simeq\widetilde{{\C}^2} $
 of the exceptional divisor which have fixed charge $i$ (cf. definition 2.1).
 We define $${\mathcal 
N}_i^f: = \{framed \,\,bundles\,\, on\,\, N(\ell) \,\, having 
\,\,charge\,\, i\}/\sim.$$

King [Kn] showed that ${\mathcal N}^f_i$ can be identified with the moduli 
space $M_i(\Sigma_1)$ of bundles on the first Hirzebruch surface 
$\Sigma_1= {\mathbb P({\mathcal O}(1)}\oplus {\mathcal O})$ having second Chern 
class $i,$ and that are trivial (and framed) over the line at infinity 
$\ell_{\infty}.$ Using this identification, we show that $M_i(\Sigma_1)$ 
is smooth and finite dimensional.

\noindent This section contains a series of lemmas, which prove the 
following result.

\vspace{3mm}

\noindent {\sc Theorem 6.1} {\it ${\mathcal N}_i^f$ is a smooth complex 
manifold.}

\vspace{3mm}

\noindent {\em 6.1 Dimension of local moduli}

%The exact dimensions of the local moduli spaces are not necessary for our 
%proof, it is only the codimensions of strata that are necessary in the 
%proof. Nevertheless, we include some examples of dimension calculation for 
%completeness.
Note that bundles on the neighborhood $N(\ell)$ of the 
exceptional divisor are completely determined by a finite infinitesimal 
neighborhood of $\ell,$ consequently the extension class and the 
inequivalent reframings depend only on a finite number of parameters, 
hence framed local moduli are finite dimensional. 

\vspace{3mm}

\noindent{\sc Proposition 6.2} {\it
Let $F:=(j,p)$ be a bundle on $N(\ell)$ with splitting type $j$ 
and extensions class $p,$ together with a trivialisation on $N^0.$
Suppose  $m>0$ is 
the $u$-multiplicity of $p$ (that is, the largest power of $u$ 
that divides $p)$. Then the dimension of the local moduli 
space at $F$ is $m (2j-(m+1)/2)$}.
\vspace{3mm}

\noindent {\it Proof.}
{ Following  [Kn], we identify ${\mathcal N}_{\bf c}^f $ 
with the moduli space $M_{\bf c}(\Sigma)$ of bundles on the 
first Hirzebruch surface with $c_1=0$ and $c_2={\bf c}$
 with a fixed trivialisation
at the line at infinity $\ell_{\infty}$. 
 By [L, Thm. 4.6], the Zariski tangent space of 
$M_{\bf c}(\Sigma)$ at $E$ is $H^1(\Sigma, \mbox{End} E \otimes 
{\mathcal O}_{\Sigma}(-\ell_{\infty})).$

We claim that  is $H^1(\Sigma, \mbox{End} E \otimes 
{\mathcal O}_{\Sigma}(-\ell_{\infty}))= H^1(N(\ell), End F)$ 
where $F = E\vert_{N(\ell)}.$ Clearly 
$H^1(N(\ell) , \End F \otimes 
{\mathcal O}_{\Sigma}(-\ell_{\infty})) = 
H^1(N(\ell), \End F)$ because ${\mathcal O}_{\Sigma}(- \ell_{\infty}) $
is trivial over $N(\ell).$ We now write 
$\Sigma = N(\ell) \cup_{N^0}N(\ell_{\infty})$ and set 
${\mathcal G} = \End E \otimes {\mathcal O}(-\ell_{\infty}).$
By Mayer-Vietoris
$$H^0(N(\ell), {\mathcal G}) \oplus H^0(N(\ell_{\infty}), {\mathcal G})
\rightarrow H^0(N^0,{\mathcal G}) \rightarrow H^1(\Sigma, {\mathcal G})
\rightarrow  $$
$$\rightarrow 
H^1(N(\ell), {\mathcal G}) \oplus H^1(N(\ell_{\infty}), {\mathcal G})
\rightarrow H^1(N^0,{\mathcal G}).$$
Here ${\mathcal G}\vert_{N(\ell_{\infty})} = {\mathcal O}^{\oplus 4} \otimes 
{\mathcal O}(-\ell_{\infty})$ and ${\mathcal G}\vert_{N^0}= 
{\mathcal O}^{\oplus 4}.$  
It then follows that the map 
$H^1(N(\ell_{\infty}), {\mathcal G}) \rightarrow H^1(N^0,{\mathcal G}) $ 
is an isomorphism, and also the map $H^0(N(\ell),{\mathcal G}) 
\rightarrow H^0(N^0,{\mathcal G})$ is an isomorphism.
The Mayer-Vietoris sequence
becomes 
$$0 \rightarrow H^1(\Sigma, {\mathcal G})
\rightarrow  
H^1(N(\ell), {\mathcal G}) 
\rightarrow 0.$$
It remains to calculate $H^1((N(\ell),{\mathcal G}) = H^1(\End F).$
The transition matrix for $\mbox{End} F $ in our canonical coordinates 
is 
$$T=
 \left( \begin {matrix} z^{2j} & - pz^j & pz^j & -p^2 \cr 
                         0 & 1 & 0 & pz^{-j} \cr
                         0 & 0 & 1 & -pz^{-j} \cr
                         0 & 0 & 0 & z^{-2j}  \end{matrix} \right)
.$$
Because ${\mathcal O}_{\Sigma}(-\ell_{\infty})$ is trivial on $N(\ell)$ 
 the transition matrix 
for $\mbox{End}  F\otimes 
{\mathcal O}_{\Sigma}(-\ell_{\infty})$ is the same $T.$
Denote by $\sim$ cohomological equivalence.
To compute the $H^1$ suppose 
$\sigma = (a,b,c,d)$ is a 1-cocycle. Then it is represented in 
$U\cap V \simeq {\mathbb C}-\{0\} \times {\mathbb C}
= \{(z\neq 0, u)\}$ in the form 
$$\sigma = \sum_{i=0}^{\infty} \sum_{l= -\infty}^{\infty} \left( 
                  \begin {matrix} a \cr b \cr c \cr d \end{matrix} \right)
z^lu^i.$$
Since terms having only positive powers of $z$ are holomorphic in $U$ 
it follows that 
$$\sigma \sim \sum_{i=0}^{\infty} \sum_{l= -\infty}^{-1} \left( 
                  \begin {matrix} a \cr b \cr c \cr d \end{matrix} \right)
z^lu^i.$$
Changing coordinates, that is, calculating $T\sigma,$
 we get the following conditions for $\sigma $ to be 
a coboundary: the expressions 
\begin{enumerate} 

\item $z^{2j} a- pz^{j} b+ pz^{j} c- p^2 d$ 
\item  $            b+  pz^{-j}d $
\item   $        c- pz^{-j}d$ 
\item    $    z^{-2j}d $
\end{enumerate}
should be holomorphic in $V= \{z^{-1}, zu)\}.$ 
Expressions 2, 3 and 4 contain no positive powers of $z,$ 
therefore they are  holomorphic in $V$ and impose no 
extra conditions on $\sigma.$ The only condition is then that 
expression 1 be holomorphic in $V.$ 
Set $p=u^mp',$ then  we need to check which terms in the expression 
$$z^{2j} a+ u^mp'z^{j} b+
u^m p'z^{j} c+  u^{2m}p'^2 d$$ are holomorphic in $V,$ 
where $a,b,c,$ and $d$ are arbitrary holomorphic in $z$ and $u.$ 
Choosing $b,c,d$ appropriately we 
can remove all terms $z^lu^i$ having $i\geq m.$ 
We are let only with  
$$z^{2j} a \sim z^{2j}\sum_{i=0}^{m-1} \sum_{l= -\infty}^{-1}   a_{il}z^lu^i
= \sum_{i=0}^{m-1} \sum_{s= -\infty}^{2j-1}a_{is} z^su^i.$$
But, since terms having 
$s \leq i$ are holomorphic in $V,$ the non-zero cocycles 
come only from the terms $z^su^i$ with $s>i,$ so
$$z^{2j}a \sim \sum_{i=0}^{m-1} \sum_{s= i+1}^{2j-1}
a_{is} z^su^i. $$
Consequently, nontrivial 
cocycles are represented by sections of the form $\sigma = (a,0,0,0)$ where  
$$ a \sim \sum_{i=0}^{m-1} \sum_{s= i-2j+1}^{-1}
a_{is} z^su^i. $$
There are    $m(2j-(m+1)/2)$ nontrivial coefficients.}
\hfill\square{5}

\vspace{3mm}

\noindent {\em 6.2 Smoothness of local moduli.}

\noindent {\sc Lemma 6.3 } {\it For all $E, E'$ in ${\mathcal N}^f_i,$ the map 
$$H^0(\Hom (E,E')) \rightarrow H^0(\Hom (E\vert_{N^0}, E'\vert_{N^0})$$ is 
injective. } \vspace{3mm}

\noindent {\it Proof.} $N^0$ is open in $\widetilde{{\mathbb C}^2},$ hence 
two holomorphic functions that coincide in $N^0$ and are globally defined 
must be equal.\hfill\square{5}

\vspace{3mm} \noindent {\sc Lemma 6.4} {\it In order to study 
endomorphisms of a bundle $V \in {\mathcal N}_i$ it suffices to choose a fixed 
transition matrix $T$ for $V$ and then to consider endomorphisms fixing 
$T.$ } \vspace{3mm}

\noindent {\it Proof.} In fact, suppose there is an automorphism $\phi$ of 
$V$ taking $T$ to $\bar T.$ Then, $\phi= (X,Y)$ is given by a pair of 
transition matrices $X \in \Gamma (U)$ and $Y \in \Gamma (V)$ such that 
$YTX=\bar T.$ But, since $T$ and $\bar T$ represent the same bundle $V,$ 
there are change of coordinates $A \in \Gamma (U)$ and $B \in \Gamma (V)$ 
such that $T = B\bar TA.$ Now define a new endomorphism of $V$ by $\bar 
\phi = (XA, BY) ,$ then $\bar \phi (T) = BYTXA= B\bar T A = T.$ 
\hfill\square{5} \vspace{3mm}

\noindent We use the canonical form of transition matrix 
$T= {\tiny \left( \begin{matrix} z^j & p \cr 0 & z^{-j} \end{matrix}\right)}$
as in $(1).$
\vspace{3mm}

\noindent {\sc Lemma 6.5} {\it Let $\phi$ be
 an  automorphism of a framed bundle $(V,f) \in {\mathcal N}^f_i.$ Then 
$\phi\vert_{N^0}$ 
can be written in the form $\phi(T)= X^{-1}TX$ with $X \in \Gamma (U_0) 
\cap \Gamma (V_0)$ and $X = {\tiny
 \left( \begin {matrix} a& b \cr 0 & d \end{matrix}\right) .} $}
\vspace{1mm}

\noindent {\it Proof.} That the automorphism can be written as $T \mapsto 
X^{-1}TX$ follows simply because $V$ is trivial over $N^0.$
 Suppose $X = {\tiny \left( \begin{matrix} a& b \cr c & d \end{matrix} 
\right)}. $ From the equality $XT=TX$ we get
 $-z^{-j} c= z^jc,$ which immediately implies that $c=0.$ \hfill 
\square{5} \vspace{3mm}

\noindent {\sc Lemma 6.6} {\it Let $\phi$ be
 an automorphism of a framed bundle $(V,f) \in {\mathcal N}^f_i$ over $N(\ell) 
\simeq \widetilde{{\mathbb C}^2}.$ Then $\phi$
 can be written in the form $\phi(T) = 
X^{-1}TX$ with
  $X = {\tiny
 \left( \begin{matrix} a_0+a& b \cr 0 & d_0 + d
 \end{matrix} \right)}$ where $a,b,d \in \Gamma(U_0) \cap \Gamma(V_0)$ 
with $a,b$ and $d$ vanishing over the exceptional divisor and $a_0,d_0$ 
constants. } \vspace{3mm}

\noindent {\it Proof.} By lemma 6.5 $X= {\tiny \left( \begin{matrix} a& b 
\cr 0 & d \end{matrix} \right)}$ over $N^0$ and we must check what are the 
possible extensions of $X$ to the full coordinate charts $U$ and $V.$
Let  $ {\tiny \left( 
\begin{matrix} \bar a& \bar b \cr \bar c & \bar d \end{matrix} \right)}$
be a holomorphic extension of $X$  
to the entire $U-$ chart. We claim the extension
 is also of the form $X = {\tiny \left( \begin{matrix} \bar a& \bar b \cr 
0 & \bar d \end{matrix} \right)}.$ In fact, $ U \simeq \mathbb C^2$ and 
$\bar c= c=0$ on $U_0$ which is an open subset of $U$ hence $\bar c=0$ 
everywhere on $U.$ Similarly on $V$ we have the form $X= {\tiny \left( 
\begin{matrix} \bar \alpha & \bar \beta \cr 0 & \bar \delta
 \end{matrix} \right)}.$ Write $\bar x = x_0 (z) +x(z,u) $ for $x \in 
\{a,b,d\} $ and write $\bar y = \bar y_0(z^{-1})+y(z^{-1}, zu) $ for $y=\{ 
\alpha, \beta, \delta\}.$ Over the exceptional divisor one has $p=0=u$ and 
the equality $(*)$ becomes $$\left ( \begin{matrix} z^j a_0 & z^j b_0 \cr
                  0 & z^{-j}d_0 \end{matrix} \right ) = \left ( 
\begin{matrix} z^j \alpha_0 & z^{-j}\beta_0 \cr
                  0 & z^{-j}\delta_0 \end{matrix} \right ) $$ where 
$a_0,b_0,d_0 $ are holomorphic in $z$ whereas $\alpha_0, \beta_0, 
\delta_0$ are holomorphic in $z^{-1}.$ It immediately follows that 
$b_0=\beta_0=0$ and that $a_0=\alpha_0$ and $d_0=\delta_0$ are constants. 
\hfill\square{5} \vspace{3mm}

\noindent {\sc Lemma 6.7}\label{traceless} {\it A framed bundle $(V,f) 
\in {\mathcal N}^f_i$ has no traceless automorphisms, unless $V$ splits. } 
\vspace{3mm}

\noindent {\it Proof.} The equality $trace (\lambda X)= \lambda  trace 
(X)$ holds for any $X \in GL(2,\mathbb C)$ and constant $\lambda.$ 
Therefore, the bundle $V$ has traceless $GL(2,\mathbb C)$ automorphisms if 
and only if it has traceless $SL(2,\mathbb C)$ automorphisms. So, we 
assume $det(X)=1.$ Over the exceptional divisor $X\vert_{\ell}= {\tiny 
\left( \begin{matrix} a_0& 0 \cr 0 & d_0 \end{matrix} \right)}.$ Hence 
$a_od_0=1.$ If $trace(X)$ is zero, then we also have $a_0+d_0=0.$ It 
follows that $a_0= \pm i.$ Suppose $a_0=i,$ the other case is analogous. 
Since $X$ is traceless, it then follows that $X= {\tiny \left( 
\begin{matrix} i+a & b \cr 0 & -i-a
 \end{matrix} \right)}.$ On the other hand $det(X)= -(i+a)^2=1$ implies 
that either $a=0$ or $a=-2i.$ Therefore $X= \pm {\tiny \left( 
\begin{matrix} i & b \cr 0 & -i
 \end{matrix} \right)}.$ Now the equality $XT=TX$ gives $$\left( 
\begin{matrix} iz^j & ip+z^{-j}b \cr 0 & -iz^{-j} \end{matrix} \right)= 
\left( \begin{matrix} iz^j & -ip+z^jb \cr 0 & -iz^{-j}
 \end{matrix} \right)$$ implying $$2ip=(z^j-z^{-j})b.$$ But $b$ has only 
positive powers of $z,$ hence the r.h.s. contains
 powers of $z$ greater of equal to $j$ whereas the l.h.s. only has powers 
of $z$ strictly
 smaller than $j.$ Hence $b=0=p.$ We conclude that the only traceless 
automorphisms of $(V,f)$ are multiples of $\tiny \left( \begin{matrix} i & 
0 \cr 0 & -i. \end{matrix} \right).$ \hfill\square{5} \vspace{3mm}

\noindent {\sc Proposition 6.8} {\it $M_i(\Sigma_1)$ is smooth.} 
\vspace{3mm}

\noindent {\it Proof.} By [L] Theorem 1.1, given the conditions of lemma 
6.2 above, $M_i(\Sigma_1) $ is smooth at $E$ provided that $H^2(sl(E) 
\otimes {\mathcal O}_{\Sigma_1}(-\ell_{\infty}))=0,$ where $sl(E)$ is the 
bundle of traceless endomorphisms of $E.$ But lemma 6.7 implies that there 
are no traceless endomorphisms of $E,$ unless $E$ splits, in which case 
$sl(E)$ is trivial. The same conclusions then holds also for bundles on 
$M(\Sigma_1).$ 
By Serre duality,  
$$H^2(sl(E) \otimes {\mathcal 
O}_{\Sigma_1}(-\ell_{\infty}))=H^0(sl(E) \otimes {\mathcal 
O}_{\Sigma_1}(\ell_{\infty})\otimes K).$$
 Since 
$K \simeq {\mathcal O}(-2f-\ell - \ell_{\infty})$ where $f$ (having $f^2=0$) 
is the class of the fiber 
of $\Sigma_1$ and $\ell ^2 = -1.$ Then $  {\mathcal 
O}_{\Sigma_1}(\ell_{\infty})\otimes K \simeq {\mathcal O}(-2f-\ell),$
and consequently
$$H^0(sl(E) \otimes {\mathcal 
O}_{\Sigma_1}(\ell_{\infty})\otimes K )
= H^0(sl(E)(-2f-\ell))$$
which we claim vanishes. In fact, suppose not, then a global section 
gives  is an injection 
$$0 \rightarrow {\mathcal O} \rightarrow sl(E) \otimes
 {\mathcal O}(-2f-\ell)$$
and a corresponding short exact sequence 
$$0 \rightarrow {\mathcal O}(2f+\ell) \rightarrow sl(E) \rightarrow Q
\rightarrow 0.$$
Taking cohomology gives 
$$0 \rightarrow H^0({\mathcal O}(2f+\ell)) \rightarrow H^0(sl(E))
 \rightarrow H^0(Q) 
\rightarrow $$
But $ H^0({\mathcal O}(2f+\ell))\neq 0$ whereas 
by lemma 6.7 \ref{traceless}  $H^0(sl(E))=0, $ a contradiction.
 \hfill\square{5} 

\vspace{3mm}

\noindent {\bf 7. Proof of the conjecture }

We assume that the Atiyah--Jones conjecture holds true for $X,$ that is, 
we assume that there exist maps
 $r_k\colon {\mathfrak M}^f_k(X) \rightarrow {\mathfrak M}^f_{k+1}(X)$ 
inducing isomorphisms in homology ${r_k}_*\colon H_q({\mathfrak M}^f_k(X)) 
\rightarrow H_q({\mathfrak M}^f_{k+1}(X))$ for $ q \leq \lfloor k/2 
\rfloor -c$ where $c$ is a constant depending on $X.$ \vspace{3mm}

\noindent{\em 7.1 Statement of the conjecture}

For a 4-manifold $X,$ let ${\mathcal MI}_k(X)$ denote the moduli space of 
framed $SU(2)$ 
 instantons on $E$ with charge $k$ and let ${\mathcal B}^{\epsilon}_{k}(X)$ 
denote the space of  framed gauge equivalence classes of 
connections on $X$ with charge $k$ 
whose self-dual part of the 
curvature has norm less than $\epsilon.$ Given a point $x_0 \in X,$
by patching a small instanton on a neighborhood of $x_0,$ Taubes 
[Ta] constructed  smooth maps $$t_k^{x_0}\colon {\mathcal MI}_k(X)\rightarrow 
{\mathcal B}^{\epsilon}_{k+1}(X).$$
He also constructed  strong deformation retracts 
$$\tau_{k+1}\colon {\mathcal B}^{\epsilon}_{k+1} (X)\rightarrow {\mathcal 
MI}_{k+1}(X).$$ 
% The map $t^{x_0}_k$ takes an instanton on $X$ with 
%charge $k$ to a connection on $X$ with
% charge $k+1$ whose self-dual part of the $L^2$ norm of its curvature is 
%less than $\epsilon.$
He then  showed that  the stable limit 
$\displaystyle\lim_{k\rightarrow \infty } {\mathcal 
MI}_k$ indeed has the homotopy type of $ {\mathcal B}(X),$ the 
space of  framed gauge equivalence classes of connections on $E.$ 
Consequently, the  Atiyah--Jones conjecture in 
homology is
 equivalent to the statement that the maps $i_k = \tau_{k+1} \circ 
t_k^{x_0}$ induce isomorphisms in homology $${i_k}_*\colon H_q({\mathcal 
MI}_k(X)) \rightarrow H_q({\mathcal MI}_{k+1}(X))$$ through a range $q(k)$ 
increasing with $k.$

For a compact K\"ahler surface $Z,$ we denote by $${\mathcal K}_X\colon{\mathcal 
MI}_k(Z)\rightarrow {\mathfrak M}_k(Z)$$ the Kobayashi--Hitchin map given 
by ${\mathcal K}(\nabla =\partial +\overline \partial) = \overline{\partial},$ 
and let ${\mathcal H}_Z := {\mathcal K}_Z^{-1}$ be the inverse map. These maps induce
 real analytic isomorphisms of moduli spaces, cf. [LT]. Using this 
translation to moduli of bundles, we aim to prove that 
$H_q({\mathfrak M}^f_k(\widetilde{X})) = H_q({\mathfrak 
M}^f_{k+1}(\widetilde{X})) $ for $q \leq \lfloor k/2 \rfloor -c,$ 
assuming the 
analogous statement is true for $X.$ 
We now stratify ${\mathfrak M}^f_k(\widetilde{X})$
 and  show that the  composite
$${{\bf t}_k}\colon 
{\mathfrak M}^f_k(\widetilde{X}) \rightarrow {\mathfrak 
M}^f_{k+1}(\widetilde{X})$$
given by $${{\bf t}_k}\colon ={\mathcal 
K}_{\widetilde{X}} \circ \tau_{k+1}\circ t_k^{\widetilde{x}_0} \circ 
{\mathcal H}_{\widetilde{X}}$$ 
is homotopy equivalent to a map that 
preserves the stratifications. \vspace{3mm}

\noindent{\em 7.2 Stratifications}

We consider stratifications of the moduli spaces which induce
filtrations with an associated Leray spectral sequence. We refer to 
those as $L$-stratifications.
\vspace{3mm}

\noindent {\sc Definition} A smooth manifold $M$ is $L$-stratified if 
there is a decomposition of $M$ into disjoint submanifolds $M(K)$ such 
that\\ (1) The index set ${\mathcal K}= \{K\}$ is finite with a given fixed 
well ordering $\leq.$\\ (2) If $K_0$ is the smallest element in $({\mathcal 
K}, \leq),$ then $M(K_0)$ is an open-dense subset of $M$.\\ (3) For all $K 
\in {\mathcal K}$ the union of the submanifolds of the same or smaller order 
$$Z(K)=\cup_{K' \leq K} M(K')$$ is an open-dense submanifold of $M.$ \\ 
(4)  For all $K \in {\mathcal K}$ the normal bundle,
 $\nu(K),$ of $M(K)$ in $M$ is orientable. \vspace{3mm}

\noindent {\sc Proposition 7.2} {\it There is an $L$-stratification
 of the moduli space of framed bundles on $\widetilde{X}$ as $${\mathfrak 
M}^f_k(\widetilde{X}) \simeq \bigcup_{i=0}^{k} {\mathfrak M}^f_{k-i}(X) 
\times {\mathcal N}^f_i.\eqno(8)$$
 } \vspace{3mm}

\noindent {\it Proof.} The existence of the 
point set decomposition follows directly 
from the definition of the ${\mathcal N}_i$ together with 
proposition 4.4. The space $K_i =  {\mathfrak M}^f_{k-i}(X) 
\times {\mathcal N}^f_i,$ having constant Euler characteristic in each
factor,
is  flat as a family of bundles.
Because  the moduli of framed stable bundles ${\mathfrak 
M}^f_k(\widetilde{X})$ is fine (cf. [HL, 4.B]) the flat family 
$ {\mathfrak M}^f_{k-i}(X) 
\times {\mathcal N}^f_i$ must be obtained by pulling back the universal
bundle on ${\mathfrak 
M}^f_k(\widetilde{X})$ after possibly twisting by a line bundle.
But, comparing with the maps appearing in proposition 4.4, shows that
the twisting can chosen to be trivial.
% and the strata $K_i$ are identified 
% with subspaces of ${\mathfrak M}^f_k(\widetilde{X}).$ 
Lemma 7.3 below shows that the inclusion $K_i \hookrightarrow 
{\mathfrak M}^f_k(\widetilde{X})$ is an immersion.

 The first stratum
 $K_0={\mathfrak M}^f_k(X) \times {\mathcal N}^f_0$ equals the set of 
pull-back bundles and is dense in ${\mathfrak M}^f_k (\widetilde{X}).$ To 
show this gives an $L$-stratification we need to prove that the normal 
bundle of each stratum $K_i = {\mathfrak M}^f_{k-i}(X) \times {\mathcal 
N}^f_i$ is orientable. Using the results from sections 3.1 and 3.2, we may 
assume that ${\mathfrak M}^f_k(\widetilde{X})$ is smooth and that 
${\mathfrak M}^f_{i}(X)$ is smooth and contains only stable bundles. By 
theorem 6.1 ${\mathcal N}^f_i$ is a complex manifold. It then follows that 
each stratum $K_i$ is a complex  submanifold of ${\mathfrak 
M}^f_k(\widetilde{X})$ and therefore has orientable normal bundle. The 
remaining properties of the $L-$stratification are shown in the following 
lemma. \hfill\square{5} \vspace{3mm}

\noindent{\sc Lemma 7.3} The inclusion $K_i \hookrightarrow 
{\mathfrak M}^f_k(\widetilde{X})$ is an immersion.
\vspace{2mm}

\noindent{\it Proof.} Let $\widetilde{E} = (E,V, \phi)$ be 
image of the pair $((E,h),(V,g))$ 
% corresponding to the decomposition
%$\widetilde{X}=  X^0 \amalg_{N_0}N(\ell),$ 
where $\phi = h \circ g^{-1}$
is the composition of the framings  $h\colon N^0 \rightarrow E$
and $g\colon N^0 \rightarrow V.$
 We want to show that the map on tangent spaces
$$T_E{\mathfrak M}^f_{k-i}(X) \times T_V{\cal N}^f_i \rightarrow 
T_{\widetilde{E}} {\mathfrak M}_k^f (\widetilde{X})$$
is injective. 
Now $X= X^0 \amalg_{N(n)-\{x\}}N(x)$ and we have the exact sequences
$$0 \rightarrow \frac{H^0(N^0,\End E)}{H^0(N(x), \End E)} \rightarrow 
T_{{\mathfrak M}^f_{k-i}(X), (E,h)} \rightarrow H^1(X,\End E) $$
and
$$0 \rightarrow \frac{H^0(N^0,\End V)}{H^0(N(\ell), \End V)} \rightarrow 
T_{N(\ell), (V,g)} \rightarrow H^1(N(\ell),\End V) .$$
Note that, by construction, the moduli of pairs $((E,h), (V,g))$ 
is the same as the moduli of triples $(E,V,\phi)$ plus the moduli of framings,
$g.$ 
But, by Hartog's we have that 
$ \frac{H^0(N^0,\End E)}{H^0(N(x), \End E)} =0,$ 
from which it follows that the map on tangent spaces is injective.  
\hfill\square{5}\vspace{3mm}

\noindent For each $0 \leq n\leq k,$ set $S_n = \cup_{i\geq n} K_i$ and 
$Z_n = {\mathfrak M}^f_k(\widetilde{X}) \setminus S_n.$ \vspace{3mm}

\noindent {\sc Lemma 7.4} {\it For all $n,$ $Z_n$ is open and dense in 
${\mathfrak M}^f_k(\widetilde{X}).$ The real codimension of the stratum
    $K_i$ in ${\mathfrak M}^f_k(\widetilde{X})$ is at least $2i.$} 
\vspace{2mm}

\noindent {\it Proof.} By definition $K_i= \mathfrak M^f_{k-i}(X)\times 
{\mathcal N}^f_i.$ The first stratum is the set of pull-back bundles $K_0 = 
\{\pi^*(E), E \in {\mathfrak M}^f_k({X})\}$ is non-empty we will se 
below that it is 
 open and dense in ${\mathfrak M}^f_k(\widetilde{X}).$ It follows that 
$S_1 = K_1 \cup K_2 \cdots \cup K_k$ has (complex) codimension 1, in 
particular $K_1$ has codimension at least 1.

For $F \in {\mathfrak M}_k(\widetilde{X}),$ set ${\bf r}(F): = (R^1\pi_* 
F)_x$ and ${\bf q}(F) := \left( (\pi_*F)^{\vee\vee}/ \pi_*F\right)_x$ 
denote the stalks of the first derived image and of the quotient sheaf $Q$ 
at $x;$ and consider the sets $${\bf R}_n:=\{F\in {\mathfrak 
M}^f_k({X})\}, h^0({\bf r} (F))\geq n\}$$ and
 $${\bf Q}_n := \{F\in {\mathfrak M}^f_k({X})\}, h^0({\bf q} (F))\geq 
n\}.$$ We consider the morphism $f\colon{\mathfrak M}^f_k({X})\rightarrow 
\{x\}.$
 Since the target is just a point, any sheaf on ${\mathfrak M}^f_k({X})$ 
is $f$-flat. Since bundles in ${\mathfrak M}^f_k({X})$ are framed stable, 
by [HL, thm.\thinspace 4.B.4] we know that there exists a universal sheaf 
${\mathcal U}$ over ${\mathfrak M}^f_k({X}).$ Regarding ${\bf r }$ and ${\bf 
h}$ as functions applied to ${\mathcal U},$ it then follows from
 Grauert's semicontinuity theorem [Ha, p.\thinspace 288] that the 
functions $h^0(x,{\bf r}({\mathcal U}))$ and $h^0(x, {\bf h}({\mathcal U}))$ are 
upper-semicontinuous. Hence for each $n,$ the sets ${\bf R}_n$ and ${\bf 
Q}_n $ are closed in ${\mathfrak M}^f_k({X}).$
 
Now set ${r_0}:={\bf r},$ $q_0:={\bf q},$ and for $1 \leq n\leq k,$ set 
$r_n:= {\bf r}\vert_{{\bf R}_{n-1}}$ and
 $q_n:= {\bf q}\vert_{{\bf Q}_{n-1}}.$ Repeating the above reasoning for 
$r_1$ and $q_1,$ we get that $ S_2 = {\bf R}_1 \cup {\bf Q}_1$ (by thm 
0.2 in [BG1]) is closed in $S_1.$ It follows that the codimension of $K_2$ in 
${\mathfrak M}^f_k({X})$ is at least 2. Now use induction on $n.$ 
\hfill\square{5} \vspace{3mm}

\noindent{\em 7.3 Maps between the spectral sequences}

%Whereas it is not strictly true that Taubes map on $\widetilde{X}$ 
%preserves the stratification $(8),$ we show that the map on ${\mathfrak 
%M}_k(\widetilde{X})$ induced by Taubes map is homotopy equivalent to a map 
%that preserves the stratification. 
Here we  show that the  map
${{\bf t}_k}\colon 
{\mathfrak M}^f_k(\widetilde{X}) \rightarrow {\mathfrak 
M}^f_{k+1}(\widetilde{X})$ is homotopy equivalent to a map that preserves
the stratifications
 $\mathfrak M^f_k(\widetilde{X})= \cup 
K_i $ and $\mathfrak M^f_{k+1}(\widetilde{X})= \cup 
K'_i$ where $K_i= \mathfrak M^f_{k-i}(X)\times {\mathcal N}^f_i$ and $K'_i= 
\mathfrak M^f_{k-i+1}(X)\times {\mathcal N}^f_i $ as in $(8).$
Let $$d \colon
 \mathfrak M^f_k(\widetilde{X})\rightarrow \mathfrak M^f_k({X} 
\amalg_f 
\widetilde{{\mathbb C}^2})$$ be the map that re-writes a framed bundle 
$\widetilde{E}^f$ on $\widetilde{X}$ into its two components of the 
decomposition $d(\widetilde{E}^f)=(E^f,V^f)$ where $E^f$ is a framed 
bundle on $X$ and $V^f$ is a framed bundle on $\widetilde{{\mathbb C}^2}.$ 
\vspace{3mm}

\noindent {\sc Remark 7.4}:
By proposition 4.4, $d$ is a bijection, and the set of pairs 
$\{(E^f,V^f)\}$ is given the topology induced by this bijection. In fact, 
$d$ is just fancy way of expressing the identity map in a form that is 
convenient for our constructions; in particular, $d$ is a homeomorphism.

\vspace{3mm}
We choose a point $x_0$ in $X$ that is far from the point $x$ we blew-up,
 in the sense that $x_0 \notin \pi (N(\ell))$ and we set 
$\widetilde{x_0}=\pi^{-1}(x_0).$ Let $$t_k^{x_0}\colon 
{\mathcal MI}_k(X)\rightarrow {\mathcal B}^{\epsilon}_{k+1}(X) \,\, \mbox{and} 
\,\,t_k^{\widetilde{x_0}}\colon {\mathcal MI}_k(\widetilde{X})\rightarrow
 {\mathcal B}^{\epsilon}_{k+1}(\widetilde{X})$$ denote Taubes patching on a 
small neighborhood of $x_0$ in $X$ and of $\widetilde{x}_0$ in 
$\widetilde{X}$ respectively, and let  
$$\tau_{k+1}\colon {\mathcal B}^{\epsilon}_{k+1} (X)
\rightarrow {\mathcal MI}_{k+1} (X)\,\, \mbox{and} \,\,
\widetilde{\tau}_{k+1}\colon {\mathcal B}^{\epsilon}_{k+1} (\widetilde{X})
\rightarrow {\mathcal MI}_{k+1} (\widetilde{X})
$$ 
be the corresponding deformation retracts.
 The idea  is that by choosing $x$ to be far 
from the exceptional divisor, we keep the patching  far from the 
exceptional divisor as well. From the point of view of vector bundles, this 
implies that we can choose to increase the second Chern class in such a way 
 that it does not alter the bundle 
near the divisor. We now make this statement precise.

\vspace{3mm} \noindent{\sc Proposition 7.5} {\it The maps ${{\bf t}_k}\colon 
{\mathfrak M}^f_k(\widetilde{X}) \rightarrow {\mathfrak 
M}^f_{k+1}(\widetilde{X})$ and
 ${{\bf r}_k}\colon {\mathfrak M}^f_k(\widetilde{X}) \rightarrow 
{\mathfrak M}^f_{k+1}(\widetilde{X})$ given by $${{\bf t}_k}\colon ={\mathcal 
K}_{\widetilde{X}} \circ \tau_{k+1}\circ t_k^{\widetilde{x}_0} \circ 
{\mathcal H}_{\widetilde{X}}$$ and $${\bf r}_k \colon = d^{-1} \circ 
({\mathcal K}_X, {\mathcal K}_{\widetilde{\mathbb C}^2})
\circ \tau_{k+1}\circ (t_k^{x_0}, id)
 \circ ({\mathcal H}_X,{\mathcal H}_{\widetilde{{\mathbb C}^2}})\circ d$$ are 
homotopically equivalent. Moreover, the map ${\bf r}_k$ 
preserves the stratification.} \vspace{3mm}

\noindent {\it Proof.} The following diagram summarises the situation 
{\scriptsize $$\begin{matrix} {\mathfrak M}^f_k(\widetilde{X})& \stackrel 
{{\mathcal H}_{\widetilde{X}}}{\longrightarrow} & {\mathcal 
MI}^f_k(\widetilde{X})&
 \stackrel { t_k^{\widetilde{x}_0} }{\longrightarrow}& {\mathcal B}^{\epsilon 
}_{k+1}(\widetilde{X})&
 \stackrel { \tilde\tau_{k+1} }{\longrightarrow}& {\mathcal 
MI}^f_{k+1}(\widetilde{X})& \stackrel {{\mathcal 
K}_{\widetilde{X}}}{\longrightarrow} & {\mathfrak 
M}^f_{k+1}(\widetilde{X})&
 \cr d\downarrow & & \wr| & & \wr|& & \wr|& & d^{-1}\uparrow\cr {\mathfrak 
M}^f_k(\bar X)& \stackrel {({\mathcal H}_X,{\mathcal H}_{\widetilde{\C}^2})} 
{\longrightarrow} & {\mathcal MI}^f_k(\bar X)& \stackrel { (t_k^{x_0}, id) 
}{\longrightarrow}& {\mathcal B}^{\epsilon }_{k+1}(\bar X)&
 \stackrel {\tau_{k+1}} {\longrightarrow}& {\mathcal MI}^f_{k+1}(\bar X)& 
\stackrel {({\mathcal K}_X,{\mathcal K}_{\widetilde{\C}^2})} {\longrightarrow} & 
{\mathfrak M}^f_{k+1}(\bar X)& \cr \end{matrix}.$$} In the second row 
$\bar X$ stands for $\bar X := 
(X-\{x\})_{\amalg_{i_1=i_2}} \widetilde{{\mathbb C}^2}$ identified as 
in section 4. A bundle on $\bar X$  given by a pair $(E,V)$ of 
bundles on $X$ and $ \widetilde{{\mathbb C}^2}$ again as in section 4. 
Note that each of the vertical maps
 is a bijection. It is important to keep in mind that all framings are 
given on open sets, not just at points. The maps ${\mathcal H}, {\mathcal K}, 
\tau$ and $t$ are well known to be continuous, cf. [LT] and [Ta2].
 Note that the homeomorphism 
$d(\widetilde{E}^f) = (E^f,V^f)$ satisfies $E= 
\pi_*\widetilde{E}^{\vee\vee}$ and $V \simeq \widetilde{E}\vert_{N(\ell)} 
;$ in particular, $E\vert_{N^0} \simeq V\vert_{N^0}$ are identified via
 the framings (see Remark 7.4), that is, there are isomorphisms 
of framed bundles
$$ E^f\vert_{N^0}
\stackrel{\small i_1}{\longleftarrow} \widetilde{E}^f\vert_{N^0 } 
\stackrel{\small i_2}{\longrightarrow} V^f\vert_{N^0} .\eqno(8)$$
 Moreover, ${\mathcal H}_X\vert_{N^0}$ takes the $\bar \partial_E$ 
operator over $N^0$ to the unique $SU(2)$ connection $\nabla_E$ on 
$E\vert_{N^0}$ having this $\bar \partial_E$ as its $(0,1)$ part.  
Uniqueness is proven as in 
 [DK, Lemma \thinspace 2.1.54]. For our purposes it
is important to notice that the proof of this lemma
 is carried out on local trivialisations, uniqueness being obtained in 
each coordinate patch. Explicitly,
the matrix of 1-forms $\alpha^{\tau}_E$ giving the partial connection 
$\bar\partial_E$ is associated to the connection $\nabla_E^{\tau} =
 \alpha^{\tau}_E
-(\alpha^{\tau})_E^*.$ Similarly, 
  ${\mathcal H}_{\widetilde{\mathbb C}^2}\vert_{N^0}$ takes the 
$\bar \partial_V$ 
operator over $N^0$ to the unique $SU(2)$ Hermitian Yang--Mills
connection $\nabla_V$ on 
$V\vert_{N^0}$ having this 
$\bar \partial_V$ as its $(0,1)$ part.
The framings on bundles induce framings on the corresponding connections, and
the isomorphisms of framed 
bundles $i_1$ and $i_2$ in $(8)$ induce isomorphisms 
of framed connections
$$ \nabla_E^f\vert_{N^0}
\stackrel{\small i_1}{\longleftarrow} \nabla_{\widetilde{E}}^f\vert_{N^0 } 
\stackrel{\small i_2}{\longrightarrow} \nabla_V^f\vert_{N^0} .\eqno(9)$$
% Moreover, ${\mathcal H}_X\vert_{N^0}$ takes the $\bar \partial_E$ 
% Hence, the maps ${\mathcal H}_X$ 
%and ${\mathcal H}_{\widetilde{\mathbb C^2}}$ coincide on $N^0.$ 
It follows that $({\mathcal H}_X,{\mathcal H}_{\widetilde{\C}^2})$ is well defined 
on the moduli space (that is, on isomorphism classes)
and continuous; and we have commutativity of the first square. A similar 
reasoning shows commutativity of the last square.

In the second square, we chose $t_k^{x_0}$ to be Taubes patching of a 
nearly anti-self-dual connection
 around a point $x_0 \notin N(\ell).$ Since the patching is a local 
operation, we can chose it so that it does not change the connection on 
$N(\ell),$ that is, we may assume that $(t_k^{x_0}\nabla)\vert_{N(\ell)}= 
\nabla\vert_{N(\ell)}$ and we can frame $(t_k^{x_0}\nabla)$ accordingly 
over $N^0.$ Hence, there is an identification
 $(t_k^{x_0}\nabla) \simeq \nabla$ over $N^0$ via the framings
and it follows that 
$(t_k^{x_0},id) $ is well defined and continuous. Using the isomorphisms 
given by the vertical maps, we may consider the maps $t:=(t_k^{x_0},id) $ 
and
 $\tilde{t}:=t_k^{\widetilde{x}_0}$ as being both defined on the same 
spaces
 ${\mathcal MI}^f_k(\widetilde{X}) \rightarrow {\mathcal B}^{\epsilon 
}_{k+1}(\widetilde{X}).$ For any instanton $\nabla \in {\mathcal 
MI}_k(\widetilde{X}),$   the 
self-dual parts of the curvatures of 
$\tilde{t}(\nabla)$ and  $t(\nabla)$ have 
$L^2$ norm  less than $\epsilon.$
 Note also that $t$ and $\tilde t$ can be chosen to change 
the connection only on a small ball $U(x_0)$ around $x_0.$ We then use the 
norm of the curvature to estimate the norm of the connections over 
$U(x_0),$ as in [DK, \S 2.3]. This shows that $t(\nabla)$ and $\tilde 
t(\nabla)$ are at distance less than $2 \epsilon$ on $U(x_0)$ and by 
construction, they coincide outside $U(x_0).$ It follows that the $L^2$ 
distance between $t$ and $\tilde{t}$ is less than $2\epsilon.$ Choosing 
$\epsilon$ small enough, this implies that $t$ and $\tilde{t}$ are 
homotopic.

In the third square of our diagram, the maps 
${\tau}_{k+1}$ and $\tilde{\tau}_{k+1}$ are deformation retractions.
Note that the map ${\tau}_{k+1}$
is applied to connections of the form $\nabla \in im(t_k^{x_0},id),$
where  $\nabla$ is anti-self-dual over $N(\ell) \simeq \widetilde{{\mathbb 
C}^2},$ and framed over $N^0.$ 
The point here is to guaranty that the image consists of 
connections that are also framed on $N^0,$
so the image falls in  ${\mathcal MI}^f_{k+1}(\bar{X}).$
In fact, this follows from the definition of  ${\tau}_{k+1}$  because 
this deformation retraction is constructed via
an application of the contraction mapping 
theorem; consequently, $im{\tau}_{k+1}$ consists of instantons that are
also  framed on $N^0.$ There is no reason to 
expect that the vertical isomorphisms 
make the diagram commutative, because there is no obvious  
way to compare  $im{\tilde\tau}_{k+1}$ to 
$im{\tau}_{k+1}.$ \
However,  
% Therefore the composites $t_k^{\tilde x_0} \circ 
%{\tau}_{k+1}$ and $t_k^{x_0}\circ {\tilde\tau}_{k+1}$ are well defined as 
%maps ${\mathcal MI}^f_{k}(\widetilde{X}) {\rightarrow} {\mathcal 
%MI}^f_{k+1}(\widetilde{X}).$ 
the horizontal maps are deformation retractions, and the vertical maps 
are isomorphisms, so the third square homotopy commutes; and this is all 
we need.
\hfill\square{5} \vspace{3mm}

\noindent {\em 7.6 Proof of AJ for blow-ups}

 Let ${\mathfrak S}_k $ denote 
the Leray spectral sequence associated to the stratification $\{K_i\}$ of 
${\mathfrak M}^f_k(\widetilde{X})$ and let ${\mathfrak S'}_{k+1} $ denote 
the Leray spectral sequence associated to the stratification
 $\{K'_i\}$ of ${\mathfrak M}^f_{k+1}(\widetilde{X})$. \vspace{3mm}

\noindent {\sc Theorem 7.6} {\it The maps ${{\bf t}_k}\colon
 {\mathfrak M}^f_k(\widetilde{X}) \rightarrow {\mathfrak 
M}^f_{k+1}(\widetilde{X})$ induce isomorphisms in homology $H_q({\mathfrak 
M}^f_k(\widetilde{X})) \rightarrow H_q({\mathfrak 
M}^f_{k+1}(\widetilde{X}))$ for $q \leq \lfloor k/2 \rfloor -c.$ } 
\vspace{3mm}

\noindent {\it Proof.} By proposition 7.5 we may assume that
 ${{\bf t}_k}\colon
 {\mathfrak M}^f_k(\widetilde{X}) \rightarrow {\mathfrak 
M}^f_{k+1}(\widetilde{X})$ respects the stratifications and therefore induces 
a map of spectral sequences
 ${{\bf t}_k}\colon {\mathfrak S}_k \rightarrow {\mathfrak S'}_{k+1}. $ 
The $E^1$ term of $ {\mathfrak S}_k $ is $E^1_{p,q} = H_p(T \nu(K_q)),$
where $T\nu$ denotes the Thom space of the normal bundle, 
and the $E^1$ term of
 $ {\mathfrak S'}_{k+1} $ is $E^1_{p,q} = H_p(T \nu(K'_q)).$ By 
hypothesis Atiyah--Jones holds for $X,$ hence we are assuming that for 
each $k$ and for $q \leq \lfloor k/2 \rfloor -c,$ $ H_q(\mathfrak 
M^f_k(X))=H_q(\mathfrak M^f_{k+1}(X)).$ Consequently $$H_q(K_i)=H_q( 
\mathfrak M^f_{k-i}(X)\times {\mathcal N}^f_i) =
 H_q( \mathfrak M^f_{k-i+1}(X)\times {\mathcal N}^f_i) =H_q(K'_i)$$ for $q\leq 
\lfloor (k-i)/2 \rfloor -c.$
 The Thom isomorphism gives $H_*(T \nu(K_i))= H_{*-\tau (K_i)} (K_i), $ 
where $\tau (K_i)$ is the real 
codimension of $K_i$ and similarly for $K'_i.$ By lemma 7.4 the real 
codimension of $K_i$ at least $2i.$ Therefore $$H_q(T \nu(K_i))=H_q(T 
\nu(K'_i))\, \mbox{ for} \, q \leq \lfloor (k-i)/2 \rfloor -c + 2 i. $$ 
The $E^1$ terms are:

$$\begin{array} {|c|c|c|c}
\cdots & \cdots& \cdots & \cr
\hline \cr
H_2(T_{\nu}K_0) & H_2(T_{\nu}(K_1) & H_2(T_{\nu}(K_2)& \cdots \cr
\hline \cr
H_1(T_{\nu}K_0) & H_1(T_{\nu}(K_1) & H_1(T_{\nu}(K_2)& \cdots \cr
\hline \cr
H_0(T_{\nu}K_0) & H_0(T_{\nu}(K_1) & H_0(T_{\nu}(K_2)& \cdots \cr
\hline 
\end{array}$$

Since $ \lfloor (k-i)/2 \rfloor -c + 2 i \geq \lfloor k/2 \rfloor -c, $ we 
conclude that
 ${{\bf t}_k} \colon {\mathfrak S}_k \rightarrow {\mathfrak S}'_{k+1} $ 
induces an isomorphism of $E^1$ terms for all $E^1_{r,s}$ with $r+s \leq 
\lfloor k/2 \rfloor -c.$ \hfill\square{5}
 
\vspace{5mm} \noindent {\sc Corollary 7.7} {\it The Atiyah--Jones 
conjecture is true for rational surfaces.} \vspace{3mm}

\noindent {\it Proof.} Every rational surface is obtained by blowing up 
points on ${\mathbb P}^2$ or on a rational ruled surface, and on these 
cases the conjecture holds true by [BHMM] and [HM].\hfill \square{5} 
\vspace{5mm}

\centerline{\large \bf References} \vspace{5mm}

\noindent [AJ] Atiyah, M. F. -- Jones, J. D. S. {\em Topological aspects 
of Yang--Mills theory}. Comm. Math. Phys. {\bf 61} (1978) no. 2, 
97--118.\\

\noindent [BG1] Ballico, E.-- Gasparim, E. {\em Vector Bundles on a 
Neighborhood of an Exceptional Curve and Elementary Transformations},
 Forum Math. {\bf 15} no.\thinspace 1, 115--122\\
 
\noindent [BG2] Ballico, E.-- Gasparim, E. {\em Numerical Invariants for 
Bundles on Blow-ups}, Proc. of the Amer. Math. Soc.
 {\bf 130} (2002) no.\thinspace 1, 23--32\\

\noindent [BHMM] Boyer, C. P.-- Hurtubise, J. C.-- Mann, B. M.-- Milgram, 
R. J. {\em The topology of instanton moduli spaces. I. The Atiyah--Jones 
conjecture}, Ann. of Math. (2) {\bf 137} (1993) no.\thinspace 3, 
561--609\\

\noindent [Do] Donaldson, S. K. {\em Polynomial invariants for smooth 
four-manifolds}, Topology {\bf 29} (1990) no.\thinspace 3, 257--315\\

\noindent [DK] Donaldson, S. K. -- Kronheimer, P. B. {\em The geometry of 
four-manifolds}, Oxford University Press, Oxford (1990)\\

\noindent [FM] Friedman, R.-- Morgan, J. {\em Certain Algebraic Surfaces 
II},
 Journal of Diff. Geometry {\bf 27} (1988) 297-369 \\

\noindent [G1] Gasparim, E. {\em Chern Classes of Bundles on Blown-up 
Surfaces}, Comm. Algebra {\bf 28} (2000) no.\thinspace 10, 4919--4926\\

\noindent [G2] Gasparim, E. {\em Holomorphic bundles on ${\mathcal O}(-k)$
 are algebraic}, Comm.  Algebra {\bf 25} (1997) no.\thinspace 9, 
3001--3009\\

\noindent [G3] Gasparim, E. {\em Holomorphic Vector Bundles on the Blow-up 
of ${\bf C}^2,$} Journal of Algebra, {\bf 199} (1998) 581--590 \\

\noindent [G4] Gasparim, E. {\em On the Topology of Holomorphic Bundles},
 Bol. Soc. Parana. Mat. (2) {\bf 18} (1998) no.\thinspace 1-2,
 113--119 \\

\noindent [G5] Gasparim, E. -- Ontaneda, P. {\em Three
applications of instanton numbers}, Comm. Math. Phys. {\bf 270} 
(2007) no. 1: 1-12 \\

\noindent [GS] Gasparim, E. -- Swanson, I. {\em Computing instanton 
numbers of curve singularities}, J. Symbolic Comput.{\bf 40}
 (2005),  no. 2, 965--978\\

\noindent{[Ha]} Hartshorne, R. {\it Algebraic Geometry.
 Graduate Texts in Mathematics} {\bf 56} Springer Verlag (1977)\\

\noindent [Hu] Hurtubise, J. {\em Instantons and jumping lines}, Comm. 
Math. Phys. {\bf 105} (1986) no. \thinspace 1, 107--122\\

\noindent [HM] Hurtubise, J. -- Milgram, R. J. {\em The Atiyah--Jones 
conjecture for ruled surfaces}, J. Reine Angew. Math. {\bf 466} (1995) 
111--143 \\

\noindent [HL] Huybrecht, D. -- Lehn, M. {\em The geometry of moduli 
spaces of sheaves}, Aspects of Mathematics, E31. Friedr. Vieweg \& Sohn, 
Braunschweig, 1997.  }\\

\noindent [Kn] King, A. {\em Instantons and holomorphic bundles on the 
blown-up plane}, Ph.D. Thesis, Oxford (1989)  \\

\noindent [Ki] Kirwan, F. {\em On spaces of maps from Riemann surfaces to 
Grassmanians and applications to the cohomology of vector bundles}, Ark. 
Math. {\bf 24} (1986) no.\thinspace 2, 221--275\\

\noindent [L] L\"ubke, M.  {\em The analytic moduli space of framed 
vector bundles}, J. reine angew. Math. {\bf 441} (1993) 45--59\\

\noindent [LT] L\"ubke, M. -- Teleman, A. {\em The Kobayashi--Hitchin 
correspondence}, World Scientific Publishing Co., Inc.,
 River Edge, NJ (1997) \\

\noindent [M] Maruyama. M. {\em Moduli of Stable Sheaves. I.} J. Math. 
Kyoto Univ. {\bf 17} (1997), no.\thinspace 1, 91--126\\

\noindent [Ta1] Taubes, C. H. {\em Path-connected Yang--Mills moduli 
spaces}, J. Diff. Geom. {\bf 19} (1984) 337--392\\

\noindent [Ta2] Taubes, C. H. {\em Metrics, connections and gluing 
theorems}, CBMS {\bf 89}, Conference on Analytic Gauge Theory, New Mexico 
State University (1994)\\

\noindent Elizabeth Gasparim,

\noindent New Mexico State University, and
 The University of Edinburgh

\end{document}